\documentclass{amsart}
\usepackage{amsmath,amssymb,latexsym}

\newcommand{\bZ}{\mbox{${\mathbb Z}$}}

\newcommand{\bC}{\mbox{${\mathbb C}$}}

\newcommand{\case}[3]{ \left\{ \begin{array}{ll} #1 &\mbox{if $#2$} \\#3 &\mbox{otherwise\,.} \end{array} \right.}
\newcommand{\casefour}[4]{ \left\{ \begin{array}{ll} #1 &\mbox{and $\;\;\,#2\,,\;\;$ or} \\#3 &\mbox{and $\;\;\,#4$\,.} \end{array} \right.}
\newcommand{\casethree}[6]{ \left\{ \begin{array}{ll} #1 &\mbox{if}\spa #2 \\ [.05in] #3 &\mbox{if}\spa #4 \\ [.05in] #5 &\mbox{if}\spa #6\,. \end{array} \right.}

\newcommand{\stacksum}[2]{\sum_{\begin{array}{c}\vspace{-6mm}\;\\ \vspace{-1mm}\scriptstyle{#1}\\ \scriptstyle{#2}\end{array}} }
\newcommand{\stackprod}[2]{\prod_{\begin{array}{c}\vspace{-6mm}\;\\ \vspace{-2mm}\scriptstyle{#1}\\ \scriptstyle{#2}\end{array}} }

\newcommand{\spa}{\;\:}

\newcommand{\fb}{\overline{F}}
\newcommand{\eb}{\overline{E}}
\newcommand{\ft}{\widetilde{F}}
\newcommand{\eh}{\widehat{E}}
\newcommand{\gb}{\overline{G}}

\newcommand{\fg}{{\mathfrak G}}
\newcommand{\fs}{{\mathfrak S}}
\newcommand{\eg}{{\rm end}(\gamma)}
\newcommand{\e}{{\mathcal{E}}}

\newcommand{\ZZ}{\mathbb{Z}}
\newcommand{\fh}{{\mathfrak H}}

\newtheorem{conj}[equation]{Conjecture}
\newtheorem{prop}[equation]{Proposition}
\newtheorem{dfn}[equation]{Definition}
\newtheorem{alg}[equation]{Algorithm}
\newtheorem{thm}[equation]{Theorem}
\newtheorem{cor}[equation]{Corollary}
\newtheorem{lem}[equation]{Lemma}
\newtheorem{exa}[equation]{Example}
\newtheorem{exas}[equation]{Examples}

\theoremstyle{remark}
\newtheorem{rem}[equation]{Remarks}
\newtheorem{rema}[equation]{Remark}

\makeatletter
\let\choose\@@choose

\makeatother

\makeatletter\@addtoreset{equation}{section}\makeatother

\begin{document}
\bibliographystyle{plain}
\title{Quantum Grothendieck Polynomials}
\author{Cristian Lenart}
\address{Department of Mathematics and Statistics, State University of New York at Albany, Albany, NY 12222, USA}
\email{lenart@albany.edu}
\urladdr{http://math.albany.edu/math/pers/lenart}
\author{Toshiaki Maeno}
\address{Department of Electrical Engineering, Kyoto University, Sakyo-ku, Kyoto 606-8501, Japan}
\email{maeno@kuee.kyoto-u.ac.jp}

\thanks{Cristian Lenart was supported by National Science Foundation grant DMS-0403029.\\ 
\indent Toshiaki Maeno was supported by Grant-in-Aid for Scientific Research.}


\begin{abstract} 
Quantum $K$-theory is a $K$-theoretic version of quantum cohomology, which was recently defined by Y.-P. Lee. Based on a presentation for the quantum $K$-theory of the classical flag variety $Fl_n$, we define and study quantum Grothendieck polynomials.  We conjecture that they represent Schubert classes (i.e., the natural basis elements) in the quantum $K$-theory of $Fl_n$, and present strong evidence for this conjecture. We describe an efficient algorithm which, if the conjecture is true, computes the quantum $K$-invariants of Gromov-Witten type for $Fl_n$. Two explicit constructions for quantum Grothendieck polynomials are presented. The natural generalizations of several properties of Grothendieck polynomials and of the quantum Schubert polynomials due to Fomin, Gelfand, and Postnikov are proved for our quantum Grothendieck polynomials. For instance, we use a quantization map satisfying a factorization property similar to the cohomology quantization map, and we derive a Monk-type multiplication formula. We also define quantum double Grothendieck polynomials and derive a Cauchy identity. Our constructions are considerably more complex than those for quantum Schubert polynomials. In particular, a crucial ingredient in our work is the Pieri formula for Grothendieck polynomials due to the first author and Sottile. 
\end{abstract}

\maketitle

\section{Introduction}

Classically, Schubert calculus is concerned with enumerative problems in
geometry, such as counting the lines satisfying some generic
intersection conditions. 
This enumeration is accomplished via a calculation in the
cohomology ring of the space of potential solutions, such as a
Grassmannian. 
The cohomology ring of a Grassmannian is well-understood combinatorially.
Less understood, particularly in combinatorial terms, are extensions to more
general flag varieties and to more general cohomology theories, such as
equivariant cohomology, quantum cohomology, or $K$-theory.
The ``modern Schubert calculus'' is concerned with the geometry and
combinatorics of these extensions. 

In this paper, we will be concerned with the variety $Fl_n$ of complete flags in ${\mathbb C}^n$. This variety (like other flag varieties) has an algebraic Schubert cell
decomposition.
Consequently, the cohomology classes and the classes of structure sheaves of {Schubert varieties} (for short, {\em Schubert
classes}) form an integral basis of the cohomology and Grothendieck rings of $Fl_n$, respectively; these classes are 
indexed by permutations in the symmetric group $S_n$. Lascoux and Sch\"utzenberger defined polynomial representatives for Schubert classes in cohomology \cite{lasps} and $K$-theory \cite{lassfm}, respectively. These polynomials, called {\em Schubert} and {\em Grothendieck polynomials}, were studied extensively; remarkable algebraic and combinatorial properties of them were discovered (see \cite{fulyt,lasagv,lrsgpp,macnsp,manfsp}). 

Motivated by ideas from string theory, mathematicians defined, for any K\"ahler algebraic manifold $X$, the (small) {\em quantum cohomology ring} $QH^*(X,\ZZ)=QH^*(X)$, which is a certain deformation of the classical cohomology ring. Fomin, Gelfand, and Postnikov \cite{fgpqsp} defined {\em quantum Schubert polynomials} using a purely algebraic and combinatorial framework. More precisely, they defined a {\em quantization map} with a nice factorization property (see (\ref{defqmapcoh})), and defined quantum Schubert polynomials as the images of the classical Schubert polynomials under this map. The quantum Schubert polynomials specialize to the classical ones upon setting the deformation parameters to 0. Several properties of quantum Schubert polynomials were derived, such as a {\em Monk-type multiplication formula} (see Theorem \ref{monkschub}). Furthermore, based on their work, as well as on a piece of geometric information in \cite{fonqcf}, Fomin, Gelfand, and Postnikov also showed that the quantum Schubert polynomials represent Schubert classes in $QH^*(Fl_n)$. {\em Quantum double Schubert polynomials} were defined and studied in \cite{cafqds,fulusp,kamqds}. 

In the recent paper \cite{leeqkt}, Y.-P. Lee defined the (small) {\em quantum $K$-theory} of a smooth projective variety $X$, denoted by $QK(X)$. This is a deformation of the ordinary $K$-ring of $X$, analogous to the relation between quantum cohomology and ordinary cohomology. The deformed product is defined in terms of certain generalizations of {\em Gromov-Witten invariants}, called {\em quantum $K$-invariants of Gromov-Witten type}. The flag variety $Fl_n$ was the first variety for which the quantum $K$-theory was studied. Givental and Lee \cite{galqkt} made the first step in computing $QK(Fl_n)$, and the complete presentation of this ring was found by Kirillov and the second author \cite{kamwp} (see Theorem \ref{presqk}). 

The goal of this paper is to extend the work in \cite{cafqds,fgpqsp,kamqds} on quantum Schubert polynomials to a quantum $K$-theory setting, and present potential applications to computing the quantum $K$-invariants of Gromov-Witten type for $Fl_n$. Thus, we define and study {\em quantum Grothendieck polynomials} (Definition \ref{defqgp}), which are a common generalization of Grothendieck and quantum Schubert polynomials. We use a new quantization map (Definition \ref{defkquant}, Corollary \ref{newquant}), which is based on the presentation of $QK(Fl_n)$ in \cite{kamwp}, and which has a factorization property similar to that of the cohomology quantization map. Thus, our quantum Grothendieck polynomials are different from those in \cite{kirqgp}, which were defined by applying the cohomology quantization map to Grothendieck polynomials. To give an idea about the complexity of Schubert calculus in quantum $K$-theory, let us mention that the largest Grothendieck polynomials for $n=5$ has 40 terms, the largest quantum Schubert polynomial has 57 terms, whereas the largest quantum Grothendieck polynomial has 1959 terms. Hence, the quantity of information encoded by quantum $K$-theory is much larger than the one encoded by $K$-theory and quantum cohomology combined.

We present several properties of quantum Grothendieck polynomials, which are natural common generalizations of the ones for Grothendieck and quantum Schubert polynomials. For instance, we present a Monk-type multiplication formula in terms of paths in the {\em quantum Bruhat graph} on the symmetric group (Theorem \ref{tmonk2}). This formula generalizes the one for Grothendieck polynomials in  \cite{lenktv} (i.e., the case $p=1$ of Theorem \ref{kthm}) and the one for quantum Schubert polynomials in  \cite{fgpqsp} (i.e., Theorem \ref{monkschub}). We also conjecture a more general Pieri-type multiplication formula (Conjecture \ref{qkpieri}). We then define the {\em quantum double Grothendieck polynomials} $\fg_w^q(x,y)$ (Definition \ref{qdgro}), and prove the {\em Cauchy identity} for quantum Grothendieck polynomials (Theorem \ref{cauchy}), which can be viewed as a weak version of their orthogonality. This identity  generalizes the Cauchy identity for Grothendieck polynomials due to Fomin and Kirillov \cite{fakgpy} (see also \cite[Proposition 2]{kirqgp}), as well as the Cauchy identity for quantum Schubert polynomials in \cite{cafqds,kamqds}. Furthermore, our Cauchy identity shows that, by analogy with the similar results for Grothendieck polynomials and quantum Schubert polynomials \cite{cafqds,kamqds}, the quantum Grothendieck polynomials $\fg_w^q$ in this paper can be recovered as $\fg_w^q=\fg_{w^{-1}}^q(y,x)|_{y=0}$ (Corollary \ref{d2s}). This leads to an explicit recursive construction of the quantum Grothendieck polynomials. An explicit nonrecursive construction in terms of {\em quiver coefficients} is presented in Section \ref{combform}.

We only conjecture that our quantum Grothendieck polynomials represent Schubert classes in the quantum $K$-theory of $Fl_n$ (Conjecture \ref{maincj}), because the geometry relevant to these classes has not yet been developed. In fact, the geometric difficulties related to proving this conjecture far exceed those in the quantum cohomology case, as explained in Section \ref{mainconj}. However, we present strong algebraic evidence for this conjecture (see Section \ref{mainconj}). Based on quantum Grothendieck polynomials, we also describe an efficient algorithm which, if the conjecture is true, computes the quantum $K$-invariants of Gromov-Witten type for $Fl_n$. An example of a computation based on our algorithm is presented (Example \ref{qkinvgw}). We also conjecture that the quantum $K$-invariants have alternating signs in a sense specified in Conjecture \ref{signalt}. 

As far as technical aspects are concerned, an important object in \cite{fgpqsp} is a commuting family of operators ${\mathcal X}_k$, which act on $QH^*(Fl_n)$ as Monk-type multiplication operators  (see (5.2) in \cite{fgpqsp}); these operators can be split into a ``non-quantum part'' and a ``quantum part''. The main technical difficulty of this paper is that, by contrast, no such splitting exists for the $K$-theoretic versions of the operators ${\mathcal X}_k$. Consequently, the latter have to be defined in the basis of Grothendieck polynomials, and one has to use multiplication formulas for Grothendieck polynomials in order to work with this definition. In essence, one needs to use the {\em Pieri formula} for Grothendieck polynomials in \cite{lasptf} (see Theorem \ref{kthm}). 

This paper is organized as follows. In Section \ref{prel}, we present background information on Schubert, Grothendieck, and quantum Schubert polynomials, as well as on the {\em Fomin-Kirillov quantum quadratic algebra}, which will be used later in the paper. In Section \ref{quantgroth}, we define our $K$-theoretic quantization map and our quantum Grothendieck polynomials; we also prove several basic properties of these polynomials. A combinatorial formula for them is given in Section \ref{combform}. In Section \ref{qmap}, we prove that the quantization map has an alternative description, which will be used in the next section; the proof of a technical result used in this section is postponed to Section \ref{pfmain}. In Section \ref{monkmul}, we derive our Monk-type multiplication formula. In Section \ref{mainconj}, we discuss the conjecture stating that our quantum Grothendieck polynomials represent Schubert classes in the quantum $K$-theory of $Fl_n$. In Section \ref{qdgroth}, we define {quantum double Grothendieck polynomials} and discuss the Cauchy identity.

{\em Acknowledgements.} We are grateful to Anders Buch and Yuan-Pin Lee for explaining to us some of their results that were used in this work. We also thank Alex Yong for the suggestion to find an explicit formula for the quantum Grothendieck polynomials based on the universal Grothendieck polynomials.

\section{Preliminaries}\label{prel}

All the polynomials in this paper are polynomials in variables $x_1,x_2,\ldots$, unless otherwise specified. 

\subsection{Schubert and Grothendieck polynomials}\label{sg}
Let $Fl_n$ be the variety of complete flags 
$(\{0\}=V_0\subset V_1\subset\ldots\subset V_n=\bC^n)$ in $\bC^n$. 
This algebraic manifold has dimension $\binom{n}{2}$.  
Its integral cohomology ring
$H^*(Fl_n)$ is isomorphic to $\bZ[x]/I_n$, where $\ZZ[x]:=\ZZ[x_1,\ldots,x_n]$ and the ideal $I_n$ is
generated by the nonconstant symmetric polynomials in
$x_1,\ldots,x_n$, and $x_i$ has cohomological degree 2. 
For this, the element $x_i$ is identified with the Chern class of the 
dual $L^*_i$ to the tautological line bundle $L_i:=V_i/V_{i-1}$.
The variety $Fl_n$ is a disjoint union of cells indexed by permutations $w$ in
the symmetric group $S_n$. 
The closure of the cell indexed by $w$ is the {\em Schubert variety} $X_w$, which
has codimension $\ell(w)$, the length of $w$ or the number of its inversions. 
The Schubert polynomial $\fs_w$ is
a certain polynomial representative for the cohomology class corresponding to
$X_w$. It is a homogeneous polynomial in $x_1,\ldots,x_{n-1}$ of degree $\ell(w)$
with nonnegative integer coefficients.  

The Grothendieck group $K(Fl_n)$ of complex vector bundles on $Fl_n$ 
is isomorphic to its Grothendieck group of coherent sheaves.
As abstract rings,  $K(Fl_n)$ and $H^*(Fl_n)$ are isomorphic.
Here, the variable $x_i$ is the $K$-theory Chern class $1{-}1/y_i$ of the line
bundle $L_i^*$, where $y_i$ represents $L_i$ in the Grothendieck
ring. 
The classes of the structure sheaves of Schubert varieties form a 
basis of $K(Fl_n)$. 
The class indexed by $w$ is represented by the Grothendieck polynomial
$\fg_w$. 
This inhomogeneous polynomial in $x_1,\ldots,x_{n-1}$ has lowest degree 
homogeneous component equal to the Schubert polynomial $\fs_w$. 

There are several constructions of Schubert and Grothendieck polynomials available, such as recursive constructions based on {\em  divided difference operators} (\ref{divdiff}) and {\em isobaric divided  difference operators} (\ref{isodivdiff}), nonrecursive combinatorial formulas etc. For these constructions and more details on Schubert and Grothendieck polynomials, we refer the reader to \cite{fulyt,lasagv,lrsgpp,macnsp,manfsp}.

While defined for $w\in S_n$, the Schubert and Grothendieck polynomials
$\fs_w$ and $\fg_w$ do not not depend on $n$.
Thus we may define them for $w$ in $S_{\infty}$, where 
$S_{\infty}:=\bigcup_n S_n$ under the usual inclusion 
$S_n\hookrightarrow S_{n+1}$. 
Both the Schubert polynomials $\fs_w$ and the
Grothendieck polynomials $\fg_w$ form bases of
$\bZ[x_1,x_2,\ldots]$, as $w$ ranges over $S_{\infty}$.

\begin{exas} {\rm Given $1\le p\le k<n$, consider the cycle
\[c[k,p]:=(k-p+1,k-p+2,\ldots,k+1)\,,\]
and let 
\[e_p^k=e_p(x_1,\ldots,x_k)\] be the  elementary symmetric polynomial of degree $p$ in $k$  variables. Then  
\begin{align}&\fs_{c[k,p]}=e_{p}^k\,,\\
&g_p^k:=\fg_{c[k,p]}=\sum_{i=p}^k(-1)^{i-p}\binom{i-1}{p-1}{e}_i^k\,.\label{expgpk}
\end{align}
The second formula first appeared in \cite{lencak}. 
}
\end{exas}

Consider the $\ZZ$-submodule of $\ZZ[x]$ given by
\begin{equation}\label{lnmod}L_n:=\langle x_1^{i_1} \ldots x_{n-1}^{i_{n-1}}\::\:0\leq i_j \leq n-j\rangle\,.\end{equation}

\begin{rema} It is well-known that $\fs_w,\fg_w\in L_n$  for $w\in S_n$. Furthermore, $L_n$ is a complement of the ideal $I_n$, so $L_n\simeq\ZZ[x]/I_n$ as $\ZZ$-modules under the quotient map.\end{rema}

We use the notation
\begin{equation}\label{sem}e_{p_1\ldots p_{m}}:=e_{p_1}^1\ldots e_{p_{m}}^{m}\,,\end{equation}
for $0\le p_i\le i$. These polynomials are called {\em standard elementary monomials}. The following is a standard result.

\begin{prop}\label{basesz} Each of the following form a $\ZZ$-linear basis of the module $L_n$, and their cosets form a $\ZZ$-linear basis of $\ZZ[x]/I_n$:
\begin{enumerate}
\item the standard elementary monomials $e_{p_1\ldots p_{n-1}}$;
\item the Schubert polynomials $\fs_w$ for $w\in S_n$;
\item the Grothendieck polynomials $\fg_w$ for $w\in S_n$.
\end{enumerate}
\end{prop}

The following result about transition matrices between the bases above is also standard.

\begin{prop}\label{transbases} The following transition matrices between bases for $L_n$ are triangular with $1$'s on the diagonal: 
\begin{enumerate}
\item from Grothendieck polynomials to Schubert polynomials;
\item from Schubert polynomials to the defining monomial basis in {\rm (\ref{lnmod})};
\item from standard elementary monomials to Schubert polynomials.
\end{enumerate}
\end{prop}

\begin{rem}\label{remtrmat} (1) The following is a well-known fact about the second transition matrix. The lexicographically smallest monomial in $\fs_w$ (with respect to the order $x_1>x_2>\ldots>x_n$ on the variables) is $x^{{\rm code}(w)}$ and occurs with coefficient 1. Here   
\begin{align}\label{notexp}&x^{\alpha}=x_1^{\alpha_1}\ldots x_{k}^{\alpha_{k}}\;\;\;\;\mbox{for $\alpha=(\alpha_1,\ldots,\alpha_k)$}\,,\;\;\mbox{and}\\ &{\rm code}(w)=(c_1,\ldots,c_{n-1})\,,\;\;\;\;\mbox{where $c_i=|\{j>i\::\: w_j<w_i\}|$}\,. \label{defcode}
\end{align}
 This remark leads to an efficient procedure for expanding a polynomial $F$ in the basis of Schubert polynomials. Indeed, we just iterate the following step: find the lexicographically smallest monomial $x^\alpha$ in $F$, and let $F:=F-a\fs_w$, where ${\rm code}(w)=\alpha$ and $a$ is the coefficient of $x^\alpha$ in $F$. 

(2) The last two transition matrices are closely related, as discussed in \cite{kamqds}, see also \cite[Corollary 5.5]{basssf} and the comment thereafter. 
\end{rem}



There are several known multiplication formulas for Schubert and Grothendieck polynomials. Most of them are expressed combinatorially in terms of the {\em Bruhat order} on $S_n$, which we now introduce. Let $t_{ab}$ denote the transposition of $(a,b)$. The  Bruhat order is the partial order on $S_n$ with covering relations $v\lessdot w=vt_{ab}$, where
$\ell(w)=\ell(v)+1$; we denote this by 
 \begin{equation}\label{edge}
    v\xrightarrow{\,({a,b})\,}w\,.
 \end{equation}
A permutation $v$ admits a cover $v\lessdot vt_{ab}$ with $a<b$ and
$v(a)<v(b)$ if and only if whenever $a<c<b$, then either $v(c)<v(a)$ or else
$v(b)<v(c)$. This is known as the {\it cover condition}; it is both explicitly and implicitly used several times in this paper. The $k$-Bruhat order is the suborder of the Bruhat order where the covers are restricted to those $v\lessdot vt_{ab}$ with $a\le k<b$. 

We will use the {\em Pieri formula} for Grothendieck polynomials derived in \cite{lasptf}. This is a formula for expanding the product $\fg_w\,g_p^k$ in the basis of Grothendieck polynomials, where $w$ is an arbitrary permutation. 

We recall some background on the Pieri formula for Grothendieck polynomials. We will use the following order on pairs of positive integers to compare covers
in a $k$-Bruhat order:
 \begin{equation}\label{ord}
  (a,b)\prec (c,d)\quad\mbox{if and only if}\quad  (b>d)\ \mbox{or}\ 
  (b=d\ \mbox{and}\ a< c)\,. 
 \end{equation}
This order first arose in connection to the {\em Monk formula} for Grothendieck polynomials \cite{lenktv}, which is the special case of the Pieri formula corresponding to $p=1$; in other words, the Monk formula provides the expansion of the product $\fg_w\,g_1^k$.

\begin{dfn}\cite{lasptf}\label{Pieri_chain}
 {\rm 
 A $k$-{\em Pieri chain} is a saturated chain $\gamma$ in the $k$-Bruhat order 
  \begin{equation}\label{chain}
     w=w_0\xrightarrow{\,({a_1,b_1})\,}
                   w_1\xrightarrow{\,({a_2,b_2})\,}\ 
               \cdots\ 
                      \xrightarrow{\,({a_s,b_s})\,}w_s=\eg\,,
         \;\;\;s=\ell(\gamma)\,,
  \end{equation}
 which satisfies the following two conditions.
 \begin{enumerate}
  \item[(P1)] $b_1\ge b_2\ge\dotsb\ge b_s\,$.\vspace{-2pt}
 \item[(P2)] For $i=2,\dotsc,\ell(\gamma)-1$, if $a_j=a_i$ for some $j<i$, then
       $(a_i,b_i)\prec(a_{i+1},b_{i+1})$.
 \end{enumerate}
}
\end{dfn}

For simplicity, if $w=w_0$ is known, we denote the above Pieri chain by the sequence $(a_1,b_1),\ldots,(a_s,b_s)$. We now consider Pieri chains $\gamma$ with certain covers marked, according to the rules (M1)-(M3) below. We indicate a marked cover by underlining its label:
 $w_{i-1}\xrightarrow{\,\underline{(a_i,b_i)}\,}w_i$. 
 \begin{enumerate}
  \item[(M1)] If the $i$th cover
              $w_{i-1}\xrightarrow{\,\underline{(a_i,b_i)}\,}w_i$ 
              is marked, then $a_j\ne a_i$ for $j<i\,$.
  \item[(M2)] If the $i$th cover
              $w_{i-1}\xrightarrow{\,(a_i,b_i)\,}w_i$ 
              is not marked and $i+1\le s$, then $(a_i,b_i)\prec(a_{i+1},b_{i+1})$.
  \item[(M3)] If $b_1=\dotsb=b_r$ and $a_1>\dotsb>a_r$ for some 
              $r\ge 1$, then $(a_r,b_r)$ is marked.
 \end{enumerate} 
If there are $p$ marked covers, we say that we have a $p$-marking. 

\begin{rema}\label{condp0}
 {\rm A Pieri chain always admits a $p$-marking for some $p>0$.}
\end{rema}

Given a Pieri chain $\gamma$, we denote by $m_p(\gamma)$ the integer $(-1)^{\ell(\gamma)-p}$ times the number of $p$-markings of $\gamma$. This number is always a signed binomial coefficient, cf. Corollary 1.16 in \cite{lasptf}. The usual convention related to binomial coefficients holds throughout this paper, namely $\binom{n}{k}$ is set to $0$ if $0\le k\le n$ does not hold. We can now state the Pieri formula. 

\begin{thm}\label{kthm} \cite{lasptf}
  We have that
 \begin{equation}\label{kpieri}
   \fg_w\,g_p^k=
   \sum_\gamma m_p(\gamma)\fg_{\eg}\,,
 \end{equation}
 where the sum is over all $k$-Pieri chains $\gamma$ (on the infinite symmetric group) 
 that begin at $w$.
 This formula has no cancellations.
\end{thm}

\subsection{Quantum Schubert polynomials}\label{qspfgp} 

In this section, we recall from \cite{fgpqsp} the main background information about quantum Schubert polynomials.
Let
\[\ZZ[q]:=\ZZ[q_1,\ldots,q_{n-1}]\,,\;\;\;\;\ZZ[q,x]:=\ZZ[q]\otimes\ZZ[x]\,.\]
The ring $\ZZ[q,x]$ is graded by ${\rm deg}(x_i)=1$ and ${\rm deg}(q_i)=2$. This grading is implicitly assumed, unless otherwise specified. 
Recall the module $L_n$ defined in (\ref{lnmod}) and let $L_n^q:=\ZZ[q]\otimes L_n$. The following notation will be often used: $q_{ij}:=q_iq_{i+1}\ldots q_{j-1}$. 

The {\em quantum elementary polynomials} $E_p^k$  ($0\le p\le k$) are defined via the Givental-Kim determinant, which is now introduced. Let
\[\Gamma_k:=\left(\begin{array}{ccccc}x_1&q_1&0&\cdots&0\\-1&x_2&q_2&\cdots&0\\ 0&-1&x_3&\cdots&0\\ \vdots&\vdots& \vdots&\ddots&\vdots\\0&0&0&\cdots&x_k\end{array}\right)\,.\]
The polynomial $E_p^k$ is defined as the coefficient of $\lambda^p$ in the characteristic polynomial ${\rm det}(1+\lambda \Gamma_k)$. Let $I_n^q$ be the ideal in the ring $\ZZ[q,x]$ generated by $E_1^n,\ldots,E_n^n$. 

\begin{prop}\cite{fgpqsp}\label{inqcompl}
The module $L_n^q$ is a complement of the ideal $I_n^q$, so $L_n^q\simeq\ZZ[q,x]/I_n^q$ as $\ZZ[q]$-modules under the quotient map. 
\end{prop}

\begin{thm}\label{qcfln}\cite{fonqcf,gakqcf}
The (small) quantum cohomology ring of $Fl_n$ has the following presentation:
\[QH^*(Fl_n)\simeq \ZZ[q,x]/{I}_n^q\,.\]
\end{thm}

The formula below allows us to compute the polynomials $E_p^k$ recursively:
\begin{equation}\label{rece0}
E_p^k=E_p^{k-1}+x_kE_{p-1}^{k-1}+q_{k-1}E_{p-2}^{k-2}\,.
\end{equation}

The polynomials $E_{p_1\ldots p_{m}}$ are defined as in (\ref{sem}) and are called {\em quantum standard elementary monomials}. Define a $\ZZ[q]$-linear quantization map $Q\::\:L_n^q\rightarrow L_n^q$ by
\begin{equation}\label{defqmapcoh}Q(e_{p_1\ldots p_{n-1}}):={E}_{p_1\ldots p_{n-1}}\,,\end{equation}
where $0\le p_i\le i$. 

\begin{dfn}\cite{fgpqsp} The {\em quantum Schubert polynomial} $\fs_w^q$, for $w\in S_n$, is defined by
\[\fs_w^q:=Q(\fs_w)\,\in\,L_n^q\subset\ZZ[q,x]\,.\]
\end{dfn}

We collect some basic properties of quantum Schubert polynomials in the following proposition.

\begin{prop}\label{propsc}\cite{fgpqsp} {\rm (1)} The polynomial $\fs_w^q$ is stable under the natural inclusion $S_n \hookrightarrow S_{N}$, $N>n$, that is, its definition does not change if $w\in S_n$ is regarded as an element of $S_N$. In consequence, we can define $\fs_w^q$ for $w$ in the infinite symmetric group $S_\infty$. 

{\rm (2)} The polynomial $\fs_w^q$ is homogeneous of degree $\ell(w)$. Specializing $q_1=\ldots=q_{n-1}=0$ yields the classical Schubert polynomial $\fs_w$.

{\rm (3)} Each of the following form a $\ZZ[q]$-linear basis of the module $L_n^q$, and their cosets form a $\ZZ[q]$-linear basis of $\ZZ[x]/I_n^q$:
\begin{itemize}
\item the quantum standard elementary monomials $E_{p_1\ldots p_{n-1}}$;
\item the quantum Schubert polynomials $\fs_w^q$ for $w\in S_n$.
\end{itemize}
\end{prop}

An algorithm for finding the expansion of an element $F\in\ZZ[q,x]/I_n^q$ in the basis of cosets of quantum Schubert polynomials was given in \cite[Corollary 12.4]{fgpqsp}. This algorithm is based on the orthogonality of quantum Schubert polynomials, which is proved in \cite{fgpqsp}, and on Gr\"obner bases techniques; it works by examining all permutations $w\in S_n$ and by finding the coefficient corresponding to the coset of $\fs_w^q$ in the expansion of $F$. Here we present a straighforward algorithm for expanding a polynomial $F\in \ZZ[q,x]$ in the basis of quantum Schubert polynomials. Our algorithm easily leads to a proof of Proposition \ref{propsc} (3) which is different from the one in \cite{fgpqsp}, where a straightening procedure is used. Without loss of generality, we will assume that the polynomial $F$ is homogeneous.

\begin{alg}\label{expsch}\hfill\\
\indent {\rm Step 1:} Let $L:=\emptyset$.\\
\indent {\rm Step 2:} Let $F_0$ be the polynomial containing all the monomials in $F$ of lowest degree with respect to the $q$ variables (with the same coefficients as in $F$). Write $F_0=q^{d_1}F_1+\ldots +q^{d_k}F_k$, using the notation in {\rm (\ref{notexp})}, where $F_i\in\ZZ[x]$. \\
\indent {\rm Step 3:} Find the expansions $F_i=\sum_{j=1}^{m_i} c_{ij}\fs_{w_{ij}}$, for $i=1,\ldots,k$, by the algorithm described in {\rm Remark  \ref{remtrmat} (1)}.\\
\indent {\rm Step 4:} Let $L:=L,(c_{11}q^{d_1},w_{11}),\ldots,(c_{km_k}q^{d_k},w_{km_k})$.\\
\indent {\rm Step 5:} Let $F:=F-\sum_{i=1}^k q^{d_i}\left(\sum_{j=1}^{m_i} c_{ij}\fs_{w_{ij}}^q\right)$.\\
\indent {\rm Step 6:} If $F\ne 0$ then go to {\rm Step 2} else output the list $L$. {\rm STOP.}
\end{alg}

The list $L$ contains the information needed to expand $F$ in the basis of quantum Schubert polynomials. Note that the permutations $w_{ij}$ do not necessarily all lie in $S_n$, unless $F\in L_n^q$. The algorithm terminates because of Proposition \ref{propsc} (2); indeed, the lowest degree of a monomial in $F$ with respect to the $q$ variables strictly increases from one iteration to the next. Therefore, the list $L$ is ordered decreasingly by the lengths of the permutations $w_{ij}$. 

The geometric relevance of the quantum Schubert polynomials is given by the following theorem.

\begin{thm}\label{repschub}\cite{fgpqsp} The quantum Schubert polynomials $\fs_w^q$ are representatives for Schubert classes in  $QH^*(Fl_n)$. 
\end{thm}

There is a {\em quantum Monk formula} for $\fs_w^q\,\fs_{s_k}^q$, where $s_k$ is the adjacent transposition $t_{k,k+1}$. 

\begin{thm}\cite{fgpqsp}\label{monkschub}
For $w\in S_\infty$, we have
\[\fs_w^q\,\fs_{s_k}^q=\fs_w^q\,(x_1+\ldots+x_k)=\sum\fs_{wt_{ab}}^q+q_{cd}\fs_{wt_{cd}}^q\,,\]
where the first sum is over all transpositions $t_{ab}$ such that $a\le k<b$ and $\ell(wt_{ab})=\ell(w)+1$, and the second sum is over all transpositions $t_{cd}$ such that $c\le k<d$ and $\ell(wt_{cd})=\ell(w)-\ell(t_{cd})=\ell(w)-2(d-c)+1$.
\end{thm}

We now recall from \cite{cafqds,kamqds} the definition of the {\em quantum double Schubert polynomials} and the {\em Cauchy identity} for quantum Schubert polynomials. The divided difference operator $\partial_i$ is by definition 
\begin{equation}\label{divdiff} \partial_i =\frac{1-s_i}{x_i-x_{i+1}}\,, \end{equation} 
where $s_i$ is the transposition of the indices $i$ and $i+1.$ 
If $w\in S_n$ has a reduced decomposition $w=s_{i_1}\ldots s_{i_l},$ 
the operator $\partial_w$ is defined by $\partial_w:=\partial_{i_1}\ldots \partial_{i_l}$;  
this definition is correct, i.e., is independent of the choice of the reduced decomposition, 
because the operators $\partial_i$ satisfy the braid relations. Let $w_\circ=n,n-1,\ldots,1$ be the longest element in $S_n$, in one-line notation.

\begin{dfn}\cite{cafqds,kamqds}\label{qdschu} 
The polynomial $\fs^q_{w_\circ}(x,y)$ for the element $w_\circ \in S_n$ 
is given by the formula 
\[ \fs^q_{w_\circ}(x,y):=\prod_{i=1}^{n-1}\left( 
\sum_{j=0}^i x_{n-i}^{i-j} E_j^i(y) \right) \,, \] 
where $E_j^i(y)$ is the corresponding quantum elementary polynomial in the $y$ variables. 
For an arbitrary element $w\in S_n,$ the {\em quantum 
double Schubert polynomial} $\fs^q_w(x,y)$ is 
\[ \fs^q_w(x,y):= \partial_{w^{-1}w_\circ}^{(x)}\fs^q_{w_\circ}(x,y)\,, \] 
where the divided difference operator $\partial_{w^{-1}w_\circ}^{(x)}$ 
acts on the $x$ variables only.
\end{dfn}

\begin{rema}
Our choice of letting the divided difference operators act on the $x$ variables is different from that in \cite{cafqds,kamqds}, where these operators act on the $y$ variables. In our case, we have $\fs_w^q(x,y)|_{q=0}=\fs_w(x,y)$, where $\fs_w(x,y)$ are the {\em double Schubert polynomials} of Lascoux and Sch\"utzenberger \cite{lasccv,lasps}; the latter represent Schubert classes in the equivariant cohomology of $Fl_n$. 
\end{rema}

We now state the main results in \cite{cafqds,kamqds}, which appear in \cite{kamqds} as Theorems B and C.

\begin{thm}\cite{cafqds,kamqds}\label{schucauchy} We have
\[ \fs^q_{w_\circ}(x,y)= \sum_{w \in S_n} \fs_w \fs^q_{ww_\circ}(y)\,. \]
\end{thm}

\begin{thm}\cite{cafqds,kamqds}\label{schud2s} We have
\[\fs_w^q=\fs_{w^{-1}}^q(y,x)|_{y=0}\,.\]
\end{thm}

\subsection{The Fomin-Kirillov quantum quadratic algebra}\label{qqa}

This algebra, introduced in \cite{fakqad}, is usually defined over the polynomial ring $\ZZ[q]$, but here we prefer to define it over $R:=\ZZ[(1-q_1)^{\pm 1},\ldots ,(1-q_{n-1})^{\pm 1}]$.

\begin{dfn} The {\em quantum Fomin-Kirillov quadratic algebra} ${\e}_n^q$ (over $R$) is the  associative 
algebra defined by the  
following generators and relations: 
\begin{itemize}
\item Generators \quad $[i,j],$ $1\leq i,j \leq n,$ $i\not= j,$ 
\item Relations  \\ 
$(0)$ $[i,j]=-[j,i],$ \\ 
$(1)$ for $i<j,$ 
\[ [i,j]^2=\left\{ \begin{array}{cc} 
q_i, & \textrm{if $j=i+1,$} \\ 
0, & \textrm{if $j>i+1,$}
\end{array}
\right. \] 
$(2)$ $[i,j][k,l]=[k,l][i,j],$ if $\{i,j\}\cap\{k,l\}~=~\emptyset,$ \\
$(3)$ $[i,j][j,k]+[j,k][k,i]+[k,i][i,j]=0.$ 
\end{itemize}
\end{dfn}

Let $h_{ij}:= 1+[i,j]$ and note that it is  invertible in ${\e}_n^q$:
\[h_{ij}^{-1}=\case{\frac{1-[i,j]}{1-q_i}}{j=i+1}{1-[i,j]}\]

We define the {\em multiplicative Dunkl elements} $X_1,\ldots ,X_n \in {\e}_n^q$ 
by the formula 
\begin{equation}\label{defxk} X_k := 1-h_{k-1,k}h_{k-2,k}\ldots h_{1k}  
h_{kn}^{-1}h_{k,n-1}^{-1} \ldots h_{k,k+1}^{-1} \,. \end{equation}

It was proved in \cite{kamsna} that the multiplicative Dunkl elements in the classical quadratic algebra (which corresponds to $q_i=0$) commute. The proof consists purely of manipulations based on the Yang-Baxter equation, which is satisfied by the corresponding elements $h_{ij}$ \cite{fakqad} (once again, $q_i=0$). Since the elements $h_{ij}$ above satisfy the Yang-Baxter equation too, the quantum version of the result in \cite{kamsna} mentioned above follows.

\begin{thm}\label{comm} (cf. \cite{kamsna}) The multiplicative Dunkl elements commute. 
\end{thm}

The following result provides a realization of the quantum $K$-theory of $Fl_n$ inside the quantum quadratic algebra. Similar realizations were proved 
for the cohomology, quantum cohomology, and $K$-theory of $Fl_n$ in \cite{fakqad}, \cite{posqvp}, and \cite{kamsna}, respectively.

\begin{thm}\label{qkenq}\cite{kamwp}
The ring $QK(Fl_n)$ is isomorphic to the subring of ${\e}_n^q$ generated by $q_1,\ldots,q_{n-1}$ and $X_1,\ldots,X_n$. 
\end{thm}

Let us now mention the quantum Bruhat representation of the quantum quadratic algebra ${\e}_n^q$.

Define the action of $[i,j]$, $i<j$, on the group algebra $R[S_n]$  by 
\begin{equation}\label{actq}[i,j]\, w := \!\!\left\{ 
\begin{array}{ll}
\!\!wt_{ij}, & \!\!\textrm{if $\ell(wt_{ij})=\ell(w)+1,$} \\ 
\!\!q_{ij}wt_{ij}, & \!\!\textrm{if $\ell(wt_{ij})=\ell(w)-2(j-i)+1,$} \\ 
\!\!0, & \!\!\textrm{otherwise.}
\end{array}
\right. \end{equation}
It is verified in \cite{fakqad} that, in this way, we obtain a representation of ${\e}_n^q$ on $R[S_n]$, and thus on $R\otimes L_n$, via the map $w\mapsto \fg_w$. As an operator on $R[S_n]$, the elements $h_{ij}=1+[i,j]$  are known as {\em quantum Bruhat operators}.

The {\em quantum Bruhat graph} is the directed graph on $S_n$ with labeled edges $w\xrightarrow{\,({i,j})\,} wt_{ij}$ in the first two cases considered in (\ref{actq}). The weight $q(\pi)$ of a path $\pi$ is the product of the monomials $q_{ij}$ for all the edges in the second case. Note that the terms on the right-hand side of the Monk formula for quantum Schubert polynomials (that is, Theorem \ref{monkschub}) correspond to the neighbors of $w$ in the quantum Bruhat graph.

\section{Quantum Grothendieck polynomials}\label{quantgroth}

In order to define quantum Grothendieck polynomials, we have to define first two sets of polynomials, which are denoted by $F_p^k$ and $\widehat{E}_p^k$.

Let $0\le p\le k\le n$. Define the polynomials $F_p^k\in\ZZ[q_1,\ldots,q_n][x]$ by
$F_0^k:=1$ and, for $p\ge 1$, by
\begin{equation}\label{deffpk} F_p^k:= \stacksum{I\subseteq [k]}{|I|=p}\prod_{i\in I}(1-x_i) 
\stackprod{i\in I}{i+1\not\in I}(1-q_i)\,. \end{equation}
Whenever the condition $0\le p\le k$ is violated by the integers $p,k$, we let $F_p^k:=0$; the same convention holds for all polynomials indexed by $p,k$ which are defined below. Let $\fb_p^k:=F_p^k|_{q_k=0}$. It is useful to also define the polynomials
\begin{equation}\label{defftpk} \ft_p^k:= \stacksum{I\subseteq [k]}{|I|=p}\prod_{i\in I}x_i \stackprod{i\in I\setminus\{1\}}{i-1\not\in I}(1-q_{i-1})\,. \end{equation}
Note that $\fb_p^k$ can be obtained from $\ft_{k-p}^k=\ft_{k-p}^k(x_1,\ldots,x_k)$ by a simple substitution, as follows:
\begin{equation}\label{subst}\fb_p^k=(1-x_1)\ldots(1-x_k)\ft_{k-p}^k\left(\frac{1}{1-x_1},\ldots,\frac{1}{1-x_k}\right)\,.\end{equation} 

The formulas below allow us to compute the polynomials $\ft_p^k$, $\fb_p^k$, and $F_p^k$ recursively.

\begin{prop}\label{recf} We have the following relations:
\begin{align}
F_p^k&=\fb_p^k-q_k(1-x_k)\fb_{p-1}^{k-1}\,,\label{recf1}\\
\ft_p^k&=\ft_p^{k-1}+x_k\ft_{p-1}^{k-1}-q_{k-1}x_k\ft_{p-1}^{k-2}\,,\label{recf2}\\
\fb_p^k&=\fb_p^{k-1}+(1-x_k)\fb_{p-1}^{k-1}-q_{k-1}(1-x_{k-1})\fb_{p-1}^{k-2}\,.\label{recf3}
\end{align}
\end{prop}

\begin{proof}
In order to prove (\ref{recf1}), let us just note that
\begin{align*}
F_p^k&=\fb_p^k-q_k(1-x_k)\stacksum{k\in I\subseteq [k]}{|I|=p}\prod_{i\in I\setminus\{k\}}(1-x_i) 
\stackprod{i\in I\setminus\{k\}}{i+1\not\in I}(1-q_i)=\\
&=\fb_p^k-q_k(1-x_k)\stacksum{I\subseteq [k-1]}{|I|=p-1}\prod_{i\in I}(1-x_i) 
\stackprod{i\in I\setminus\{k-1\}}{i+1\not\in I}(1-q_i)=\fb_p^k-q_k(1-x_k)\fb_{p-1}^{k-1}\,.
\end{align*}

By splitting the sum in the right-hand side of (\ref{defftpk}) as follows, we have
\begin{align*}
\ft_p^k= \ft_p^{k-1}&+x_k\stacksum{k\in I\subseteq [k]}{|I|=p}\prod_{i\in I\setminus\{k\}}x_i \stackprod{i\in I\setminus\{1,k\}}{i-1\not\in I}(1-q_{i-1})-\\
&-q_{k-1}x_k\stacksum{k\in I\subseteq [k]\setminus\{k-1\}}{|I|=p}\prod_{i\in I\setminus\{k\}}x_i \stackprod{i\in I\setminus\{1,k\}}{i-1\not\in I}(1-q_{i-1})\,.
\end{align*}
It is easy to see that the first sum is precisely $\ft_{p-1}^{k-1}$, while the second one is $\ft_{p-1}^{k-2}$.

The recurrence relation (\ref{recf3}) follows easily from (\ref{recf2}) via the substitution (\ref{subst}). 
\end{proof}

Let 
\begin{equation}\label{defepk}\widehat{E}_p^k:=\sum_{i=0}^p(-1)^i\binom{k-i}{p-i}F_i^k\,,\;\;\;\mbox{and}\;\;\;\eb_p^k:=\eh_p^k|_{q_k=0}\,.\end{equation}
Note that, by M\"obius inversion, we also have
\begin{equation}\label{definve}F_p^k=\sum_{i=0}^p(-1)^i\binom{k-i}{p-i}\widehat{E}_i^k\,.\end{equation}
The important role played by the polynomials $\widehat{E}_p^k$ is discussed below. To be more specific, a presentation of the ring $QK(Fl_n)$ is given in terms of them, in the same way in which a presentation of $QH^*(Fl_n)$ is given in terms of the polynomials $E_i^n$ (cf. Theorem \ref{qcfln}). 

Let $\widehat{I}_n^q$ be the ideal in the ring $\ZZ[q,x]$  generated by $\eb_i^n$ for $i=1,\ldots,n$. The presentation for the ring $QK(Fl_n)$ given below is easily deduced from the one in \cite{kamwp} by a simple substitution, cf. (\ref{subst}). 

\begin{thm}\label{presqk} (cf. \cite{kamwp}) The (small) quantum $K$-theory ring of $Fl_n$ has the following presentation: 
\[QK(Fl_n)\simeq\ZZ[q,x]/\widehat{I}_n^q\,.\]
\end{thm}

The formulas below allow us to compute the polynomials $\eb_p^k$ and $\widehat{E}_p^k$ recursively.

\begin{prop}\label{rece} We have the following relations:
\begin{align}
\eh_p^k&=\eb_p^{k-1}+q_k(1-x_k)\eb_{p-1}^{k-1}\,,\label{rece1}\\
\eb_p^k&=\eb_p^{k-1}+x_k\eb_{p-1}^{k-1}+q_{k-1}(1-x_{k-1})(\eb_{p-1}^{k-2}+\eb_{p-2}^{k-2})\,.\label{rece2}
\end{align}
\end{prop}

\begin{proof} The relation (\ref{rece1}) easily follows based on (\ref{defepk}) and (\ref{recf1}). The right-hand side of (\ref{rece2}) can be rewritten as follows:

\begin{align*}
&\sum_{i=0}^p(-1)^i\binom{k-i-1}{p-i}\fb_i^{k-1}+x_k\sum_{i=0}^{p-1}(-1)^i\binom{k-i-1}{p-i-1}\fb_i^{k-1}+\\ 
+&q_{k-1}(1-x_{k-1})\left(\sum_{i=0}^{p-1}(-1)^i\binom{k-i-2}{p-i-1}\fb_i^{k-2}+\sum_{i=0}^{p-2}(-1)^i\binom{k-i-2}{p-i-2}\fb_i^{k-2}\right)=\\
=&\sum_{i=0}^p(-1)^i\binom{k-i}{p-i}\fb_i^{k-1}-\sum_{i=0}^{p-1}(-1)^i\binom{k-i-1}{p-i-1}\fb_i^{k-1}-\\
-&x_k\sum_{i=1}^p(-1)^i\binom{k-i}{p-i}\fb_{i-1}^{k-1}-q_{k-1}(1-x_{k-1})\sum_{i=1}^p(-1)^i\binom{k-i}{p-i}\fb_{i-1}^{k-2}=\\
=&\binom{k}{p}\fb_0^k+\sum_{i=1}^p(-1)^i\binom{k-i}{p-i}\fb_i^{k-1}+\sum_{i=1}^p(-1)^i\binom{k-i}{p-i}\fb_{i-1}^{k-1}-\\
-&x_k\sum_{i=1}^p(-1)^i\binom{k-i}{p-i}\fb_{i-1}^{k-1}-q_{k-1}(1-x_{k-1})\sum_{i=1}^p(-1)^i\binom{k-i}{p-i}\fb_{i-1}^{k-2}=\\
=&\sum_{i=0}^p(-1)^i\binom{k-i}{p-i}\fb_i^k=\eb_p^k\,.
\end{align*}
Here the first expression and the last equality use (\ref{defepk}), the first equality is based on applying Pascal's identity twice,  and the third equality uses (\ref{recf3}).
\end{proof}

Let us now introduce the $K$-theoretic quantization map. We define $\widehat{E}_{p_1\ldots p_{m}}$ as in (\ref{sem}). It is easy to see that $\widehat{E}_{p_1\ldots p_{n-1}}\in L_n^q$.

\begin{dfn}\label{defkquant} Let  $\widehat{Q}\::\:L_n^q\rightarrow L_n^q$ be the $\ZZ[q]$-linear map given by
\[\widehat{Q}(e_{p_1\ldots p_{n-1}}):=\widehat{E}_{p_1\ldots p_{n-1}}\,,\]
where $0\le p_i\le i$. 
\end{dfn}

We can express the quantization map in terms of the polynomials $F_p^k$ instead. Let
\begin{equation}\label{defsfpk}f_p^k:=F_p^k|_{q=0}=e_p(1-x_1,\ldots,1-x_k)\,.\end{equation}
We define $f_{p_1\ldots p_{m}}$ and $F_{p_1\ldots p_{m}}$ as in (\ref{sem}).

\begin{prop}\label{quantf}
We have
\[\widehat{Q}(f_{p_1\ldots p_{n-1}})={F}_{p_1\ldots p_{n-1}}\,,\]
where $0\le p_i\le i$. 
\end{prop}

\begin{proof}
We have 
\begin{equation}\label{esubs}f_p^k=\sum_{i=0}^p (-1)^i\binom{k-i}{p-i}e_i^k\,.\end{equation}
The definition of the quantization map and its linearity imply
\[\widehat{Q}(f_{p_1\ldots p_{n-1}})=\widehat{Q}(f_{p_1}^1)\ldots \widehat{Q}(f_{p_{n-1}}^{n-1})\,.\]
The proof is completed by noting that, based on (\ref{esubs}), Definition \ref{defkquant}, and (\ref{definve}), we have
\[\widehat{Q}(f_p^k)=\sum_{i=0}^p (-1)^i\binom{k-i}{p-i}\widehat{E}_i^k=F_p^k\,.\]
\end{proof}

We can now define our quantum Grothendieck polynomials.

\begin{dfn}\label{defqgp} The {\em quantum Grothendieck polynomial} $\fg_w^q$, for $w\in S_n$, is
\[\fg_w^q:=\widehat{Q}(\fg_w)\,\in\,\ZZ[q,x]\] 
\end{dfn}

\begin{exa}\label{exs3}{\rm 
The following is the list of the quantum Grothendieck polynomials for $S_3$. 
\begin{align*} &\fg^q_{\rm id} = 1, \;\;\;\;\;\;\; \fg^q_{213} = (1-q_1)x_1+q_1, \\ 
& \fg^q_{132}= -(1-q_2)x_1x_2 + (1-q_2)x_1 + (1-q_2)x_2 + q_2, \\
& \fg^q_{312}= (1-q_1)^2x_1^2-q_1(1-q_2)x_1x_2 - (2q_1^2+q_1q_2-3q_1)x_1 + 
q_1(1-q_2)x_2 + q_1^2+ q_1q_2-q_1, \\ 
& \fg^q_{231}= (1-q_2)x_1x_2 -(q_1-q_2)x_1+q_1, \;\;\;\;\;\;\;
\fg^q_{321}= \fg^q_{213}\fg^q_{231}. \end{align*}
}\end{exa}

We now state several basic properties of quantum Grothendieck polynomials; some of them generalize results about Grothendieck polynomials and quantum Schubert polynomials (see Sections \ref{sg} and \ref{qspfgp}). We start with a stability result that is immediate from our definitions.

\begin{prop}\label{stab} The polynomial $\fg_w^q$ is stable under the natural inclusion 
$S_n \hookrightarrow S_{N}$, $N>n$, that is, its definition does not change if $w\in S_n$ is regarded as an element of $S_N$. 
\end{prop}

It is clear from definitions that if the code of $w$ (introduced in (\ref{defcode})) is a partition $\lambda$ such that its {\em conjugate} (i.e., reflection with respect to the diagonal) $\lambda'=(\mu_1\ge\ldots\ge\mu_k>0)$ has no two parts equal, then we have 
\[\fg_w^q=\fg_{2,\ldots\mu_1+1, 1}^q\ldots \fg_{2,\ldots\mu_k+1, 1}^q\,.\]
 The following factorization property, which can be iterated in the obvious way, is also immediate from definitions. 

\begin{prop} Let $w$ be a permutation in $S_{p+q}$ such that $w=w_1w_2$, where $w_1$ and $w_2$ are the images of permutations in $S_p$ and $S_q$ under the inclusions that factor through $S_p\times\{{\rm id}\}$ and $\{{\rm id}\}\times S_q$, respectively. Then $\fg_w^q=\fg_{w_1}^q\fg_{w_2}^q$. 
\end{prop}

\begin{proof} The statement follows easily from the following classical facts:
\[\fg_w=\fg_{w_1}\fg_{w_2}\,,\;\;\;\;\fg_{w_1}\in\ZZ[x_1,\ldots,x_{p-1}]\,,\;\;\;\;\fg_{w_2}\in\ZZ[x_1,\ldots,x_p]^{S_p}\otimes\ZZ[x_{p+1},\ldots,x_{p+q}]\,.\]
\end{proof}

The Grothendieck polynomials and the quantum Schubert polynomials can both be recovered from the quantum Grothendieck polynomials. Recall the grading of the ring $\ZZ[q,x]$ defined in Section \ref{qspfgp}, namely ${\rm deg}(x_i)=1$ and ${\rm deg}(q_i)=2$. 

\begin{prop}\label{lowest} {\rm (1)} The lowest homogeneous component of $\widehat{E}_p^k$ is the quantum elementary  polynomial $E_p^k$, while its  specialization at   $q_1=\ldots=q_n=0$ is the  elementary   symmetric polynomial $e_p^k$. 

{\rm (2)} The lowest homogeneous component of the quantum Grothendieck polynomial $\fg_w^q$ is the quantum Schubert polynomial $\fs_w^q$, while its  specialization at   $q_1=\ldots=q_n=0$ is the Grothendieck polynomial $\fg_w$.
\end{prop}

\begin{proof}
The first part of the proposition follows easily from (\ref{rece1}) and (\ref{rece2}). Indeed, by taking the lowest homogeneous component on the right-hand side of (\ref{rece2}), we obtain precisely the recurrence relation (\ref{rece0}) for the polynomials $E_p^q$. A similar reasoning works when setting the variables $q_i$ to 0. The second part is immediate based on the first part and the definition of the quantum Grothendieck polynomials.
\end{proof}

Let $G_0^k:=1$ and, for $p\ge 1$, $G_p^k:=\fg_{c[k,p]}^q$. We define ${g}_{p_1\ldots p_{m}}$ and ${G}_{p_1\ldots p_{m}}$ as in (\ref{sem}). 

\begin{prop} We have 
\begin{equation}
G_p^k=\sum_{i=p}^k(-1)^{i-p}\binom{i-1}{p-1}\widehat{E}_i^k=1+\sum_{i=k-p+1}^k(-1)^{i-k+p}\binom{i-1}{k-p}F_i^k\,.\label{defgpk}\end{equation}
Given $0\le p_i\le i$, we also have
\begin{equation}\label{charq}\widehat{Q}(g_{p_1\ldots p_{n-1}})={G}_{p_1\ldots p_{n-1}}\,.\end{equation}
\end{prop}

\begin{proof}
The first equality in (\ref{defgpk}) is immediate from the definition of the quantization map and (\ref{expgpk}). The second equality, which is not used elsewhere in this paper, is left to the reader; the main idea is to express the polynomials $\widehat{E}_i^k$ in terms of $F_j^k$ based on (\ref{defepk}), which reduces the equality to a set of binomial identities. The description of the quantization map in (\ref{charq}) follows from the definition of this map, its linearity, and (\ref{defgpk}). 
\end{proof}

\begin{prop}\label{bases} The sets 
\[\{{F}_{p_1,\ldots,p_{n-1}}\::\:0\le p_i\le i\}\,,\;\;\{\fg_w^q\::\:w\in S_n\}\,,\;\;\{\widehat{E}_{p_1,\ldots,p_{n-1}}\::\:0\le p_i\le i\}\,,\;\;\{{G}_{p_1,\ldots,p_{n-1}}\::\:0\le p_i\le i\}\]
are $R$-linear bases of the module $R\otimes L_n$.
\end{prop}

\begin{proof}
It suffices to show that the first set is an $R$-linear basis of $R\otimes L_n$. Indeed, the following hold:
\begin{itemize}
\item the transition matrix from $\widehat{E}_{p_1,\ldots,p_{n-1}}$ to $F_{p_1,\ldots,p_{n-1}}$ is triangular with $\pm 1$'s on the diagonal, by (\ref{defepk});
\item the transition matrix from ${G}_{p_1,\ldots,p_{n-1}}$ to $\widehat{E}_{p_1,\ldots,p_{n-1}}$ is triangular with $1$'s on the diagonal, by (\ref{defgpk});
\item the transition matrix from $\fg_w^q$ to $\widehat{E}_{p_1,\ldots,p_{n-1}}$ is triangular with $1$'s on the diagonal, by Definition \ref{defqgp} and Proposition \ref{transbases} (1) and (3).
\end{itemize}

We will now prove that the first set is a basis. Throughout the remainder of the proof, a basis means an $R$-linear basis, and a linear combination means one with coefficients in $R$; in addition, by degree we mean the degree with respect to the $x$ variables only (i.e., the degrees of the $q$ variables are set to 0). In order to prove the claim, it suffices to show that every element in the defining monomial basis of $L_n^q$ can be expressed as a linear combination of elements ${F}_{p_1,\ldots,p_{n-1}}$. We prove this by induction on the degrees of the mentioned monomials. Let us fix a degree $k$ and expand an element ${F}_{p_1,\ldots,p_{n-1}}$ with $p_1+\ldots+p_{n-1}=k$ in terms of the $x$ and $q$ variables. Note that the highest degree component in this expansion consists precisely of the monomials in the expansion of $e_{p_1,\ldots,p_{n-1}}$, but the coefficients are now, up to sign, products of factors $1-q_i$. By Proposition \ref{transbases} (2) and (3), any monomial of degree $k$ can be expressed as a linear combination of ${F}_{p_1,\ldots,p_{n-1}}$ with $p_1+\ldots+p_{n-1}=k$ and lower degree monomials. The proof is concluded by invoking the induction hypothesis.
\end{proof}

\begin{rema}\label{winqcompl} It can be shown by an inductive argument based on Proposition \ref{inqcompl} that $R\otimes L_n$ is a complement of $\widetilde{I}_n^q:=R\otimes \widehat{I}_n^q$ in $R[x]$. However, we omit the details here. It then follows that $L_n^q\simeq R[x]/\widetilde{I}_n^q$ as $R$-modules under the quotient map. Furthermore, the cosets of the elements in the four families in Proposition \ref{bases} are $R$-linear bases of $R[x]/ \widetilde{I}_n^q$.
\end{rema}

Given a polynomial $F\in \ZZ[q,x]$ which can be written as a $\ZZ[q]$-linear combination of quantum Grothendieck polynomials (not necessarily for $S_n$), we present an efficient algorithm for finding this expansion. This algorithm does not work (i.e., does not terminate) if $F$ does not satisfy the above condition. An important geometric application of this algorithm is given in Section \ref{mainconj}, provided that the conjecture stated there is true. 

\begin{alg}\label{expgro}\hfill \\
\indent {\rm Step 1:} Let $L:=\emptyset$.\\
\indent {\rm Step 2:} Let $F_0$ be the lowest homogeneous component of $F$.\\
\indent {\rm Step 3:} Find the expansion $F_0=\sum_{i=1}^k c_i\fs_{w_i}^q$ by Algorithm {\rm \ref{expsch}}.\\
\indent {\rm Step 4:} Let $L:=L,(c_1,w_1),\ldots,(c_k,w_k)$.\\
\indent {\rm Step 5:} Let $F:=F-\sum_{i=1}^k c_i\fg_{w_i}^q$.\\
\indent {\rm Step 6:} If $F\ne 0$ then go to {\rm Step 2} else output the list $L$. {\rm STOP.}
\end{alg}

The list $L$ contains the information needed to expand $F$ in the basis of quantum Grothendieck polynomials. Note that the permutations $w_{i}$ all lie in $S_n$ precisely when $F\in L_n^q$. The algorithm terminates because of the condition on $F$ and Proposition \ref{lowest} (2); indeed, the degree of the lowest homogeneous component of $F$ strictly increases from one iteration to the next. Therefore, the list $L$ is ordered increasingly by ${\rm deg}(c_i)+\ell(w_i)$, that is, by the degrees of the lowest homogeneous components of $c_i\fg_{w_i}^q$. 

Let $\overline{G}_p^k:=G_p^k|_{q_k=0}$. The formulas below allow us to compute the polynomials $\gb_p^k$ and $G_p^k$ recursively. Upon setting $q_i=0$, these formulas specialize to Lascoux's transition formula for Grothendieck polynomials \cite{lastgp}.

\begin{prop}\label{recg} We have the following relations: 
\begin{align}
&G_p^k=\overline{G}_p^k-q_k(1-x_k)(\overline{G}_p^{k-1}-\overline{G}_{p-1}^{k-1})\,,\label{recg1}\\
&\overline{G}_p^{k}=(1-x_k)\overline{G}_p^{k-1}+x_k\overline{G}_{p-1}^{k-1}-q_{k-1}(1-x_{k-1})(\overline{G}_{p-1}^{k-2}-\overline{G}_{p-2}^{k-2})\,,\label{recg2}\\
&G_p^k-G_{p-1}^{k-1}=(1-q_k)(1-x_k)(\gb_p^{k-1}-\gb_{p-1}^{k-1})\,.\label{recg3}
\end{align}
\end{prop}

\begin{proof}
We have
\begin{align*}
G_p^k&=\gb_p^k+q_k(1-x_k)\sum_{i=p}^k(-1)^{i-p}\binom{i-1}{p-1}\eb_{i-1}^{k-1}=\\
&=\gb_p^k-q_k(1-x_k)\sum_{i=p-1}^{k-1}(-1)^{i-p}\binom{i}{p-1}\eb_{i}^{k-1}=\\
&=\gb_p^k-q_k(1-x_k)\left(\sum_{i=p}^{k-1}(-1)^{i-p}\binom{i-1}{p-1}\eb_{i}^{k-1}+\sum_{i=p-1}^{k-1}(-1)^{i-p}\binom{i-1}{p-2}\eb_{i}^{k-1}\right)=\\
&=\gb_p^k-q_k(1-x_k)(\gb_p^{k-1}-\gb_{p-1}^{k-1})\,.
\end{align*}
Here the first equality is based on (\ref{defgpk}) and (\ref{rece1}), the third one on Pascal's identity, and the last one on (\ref{defgpk}). 

The right-hand side of (\ref{recg2}) can be rewritten as follows:
\begin{align*}
&(1-x_k)\sum_{i=p}^{k-1}(-1)^{i-p}\binom{i-1}{p-1}\eb_i^{k-1}+x_k\sum_{i=p-1}^{k-1}(-1)^{i-p+1}\binom{i-1}{p-2}\eb_i^{k-1}-\\
-&q_{k-1}(1-x_{k-1})\left(\sum_{i=p-1}^{k-2}(-1)^{i-p+1}\binom{i-1}{p-2}\eb_i^{k-2}-\sum_{i=p-2}^{k-2}(-1)^{i-p}\binom{i-1}{p-3}\eb_i^{k-2}\right)=\\
=&\sum_{i=p}^{k-1}(-1)^{i-p}\binom{i-1}{p-1}\eb_i^{k-1}-x_k\sum_{i=p-1}^{k-1}(-1)^{i-p}\binom{i}{p-1}\eb_i^{k-1}+\\
+&q_{k-1}(1-x_{k-1})\left(\sum_{i=p-2}^{k-2}(-1)^{i-p}\binom{i+1}{p-1}\eb_i^{k-2}-\sum_{i=p-1}^{k-1}(-1)^{i-p}\binom{i}{p-1}\eb_i^{k-2}\right)=
\end{align*}

\begin{align*}
=&\sum_{i=p}^{k}(-1)^{i-p}\binom{i-1}{p-1}\eb_i^{k-1}+x_k\sum_{i=p}^{k}(-1)^{i-p}\binom{i-1}{p-1}\eb_{i-1}^{k-1}+\\
+&q_{k-1}(1-x_{k-1})\left(\sum_{i=p}^{k}(-1)^{i-p}\binom{i-1}{p-1}\eb_{i-2}^{k-2}+\sum_{i=p}^{k}(-1)^{i-p}\binom{i-1}{p-1}\eb_{i-1}^{k-2}\right)=\\
=&\sum_{i=p}^{k}(-1)^{i-p}\binom{i-1}{p-1}\eb_i^{k}=\gb_p^k\,.
\end{align*}
Here the first expression and the last equality use (\ref{defgpk}), the first equality is based on applying Pascal's identity twice (the second time in the form $\binom{i-1}{p-2}+\binom{i-1}{p-3}=\binom{i+1}{p-1}-\binom{i}{p-1}$),  and the third equality uses (\ref{rece2}).

The third relation (\ref{recg3}) easily follows from the previous two. Indeed, by plugging $\gb_p^k$ as given by (\ref{recg2}) into (\ref{recg1}), we obtain, after slightly rearranging the terms:
\[G_p^k=(1-q_k)(1-x_k)(\gb_p^{k-1}-\gb_{p-1}^{k-1})+\left[\gb_{p-1}^{k-1}-q_{k-1}(1-x_{k-1})(\overline{G}_{p-1}^{k-2}-\overline{G}_{p-2}^{k-2})\right]\,.\]
By (\ref{recg1}), the expression inside the square bracket is $G_{p-1}^{k-1}$.
\end{proof}

As an easy corollary of the first two relations in Proposition \ref{recg}, we compute $G_1^k=\fg_{s_k}^q$.

\begin{cor}\label{gsk}
\[G_1^k=\fg_{s_k}^q=1-(1-x_1)\ldots(1-x_k)(1-q_k)\,.\]
\end{cor}

\section{A combinatorial formula for quantum Grothendieck polynomials}\label{combform}

In this section, we present an explicit combinatorial formula for the quantum Grothendieck polynomials. The suggestion to find such a formula was made to us by A. Yong \cite{yonpc}. A similar formula for the quantum Schubert polynomials was given in \cite{bktspq}, and was based on the formula for the {\em universal Schubert polynomials} of Fulton \cite{fulusp}. By analogy, the formula presented here is based on the one for the {\em universal Grothendieck polynomials} in \cite{bkssgp,bktgpq}. 

We will now explain the background. The universal Grothendieck polynomial, denoted $\fg_w(c)$ (for $w\in S_n$), is a polynomial in independent variables $c_p(k)$. It is obtained in a similar way to a quantum Grothendieck polynomial. Indeed, we express the classical Grothendieck polynomial $\fg_w$ as a linear combination of products $g_{p_1}^1\ldots g_{p_{n-1}}^{n-1}$ (cf. Propositions \ref{basesz} (1) and \ref{transbases} (1)), and then replace each $g_p^k$ by $c_p(k)$. The following formula for $\fg_w(c)$ was given in \cite{bktgpq}:
\begin{equation}\label{quiver}
\fg_w(c)=\sum_\nu c_{w,\nu}^{(n)}\:\fg_{\nu^n}(c(n))\:\fg_{\nu^{n-1}}(c(n-1)-c(n))\:\ldots\:\fg_{\nu^1}(c(1)-c(2))\,;\end{equation}
here the sum is over finitely many sequences of partitions $\nu=(\nu^1,\ldots,\nu^n)$, and $c_{w,\nu}^{(n)}$ are special cases of {\em quiver coefficients}. A combinatorial formula for these coefficients, based on a generalization of the {\em Robinson-Schensted-Knuth insertion} algorithm (called {\em Hecke insertion}) was given in \cite{bkssgp}. 

We now define the factors on the right-hand side of (\ref{quiver}) based on two ingredients: the coproduct in the bialgebra of stable Grothendieck polynomials defined by Buch \cite{buclrr}, and Buch's {\em Jacobi-Trudi formula} for stable Grothendieck polynomials \cite{bucgcq}. We start by setting
\[\fg_\nu(c(k)-c(k-1)):=\sum_{\lambda\subseteq\nu} \fg_{\nu/\!\!/\lambda}(c(k))\:\fg_{\lambda}(-c(k-1))\,.\]
The first factor on the right-hand side is given by
\[\fg_{\nu/\!\!/\lambda}(c(k)):=\sum_{\mu\subseteq\nu}d_{\lambda\mu}^\nu\fg_{\mu}(c(k))\,;\]
here $d_{\lambda\mu}^\nu$ are certain structure constants for the multiplication of stable Grothendieck polynomials. A combinatorial formula for $d_{\lambda\mu}^\nu$ was given in \cite{buclrr} as a generalization of the classical Littlewood-Richardson rule for multiplying Schur functions. Note that $\fg_{\nu/\!\!/\lambda}(c(k))$ depends on the shapes $\nu$ and $\lambda$ themselves, not just on the skew diagram $\nu/\lambda$; in particular, $\fg_{\nu/\!\!/\nu}(c(k))=1$ if and only if $\nu$ is the empty partition. 

Finally, the polynomials $\fg_{\mu}(c(k))$ and $\fg_{\lambda}(-c(k-1))$ are computed recursively by Buch's Jacobi-Trudi formula, as explained below. Let $\mu'$ denote the conjugate of the partition $\mu$, and let $(a,\mu)$ denote the partition obtained from $\mu$ by adding a first row of length $a$. Then, by \cite[Theorem 6.1]{bucgcq}, we have
\[\fg_{(a,\mu)'}(c(k))=c_a(k)\:\fg_{\mu'}(c(k))+\stacksum{1\le q\le \min(\mu_1',k-a)}{0\le t\le  k-a-q}(-1)^q\binom{q-1+t}{t}c_{a+q+t}(k)\:\fg_{\mu'/\!\!/(q)}(c(k))\,.\]
The polynomial $\fg_{\lambda}(-c(k-1))$ is computed similarly, but now we have to replace $k$ by $k-1$ in the above formula, and $c_{b}(k)$ by $\fg_{(1^b)}(-c(k-1))$. The latter is computed recursively based on
\[\sum_{i=0}^b c_i(k-1)\,\fg_{(1^b)/\!\!/(1^i)}(-c(k-1))=\delta_{0,b}\,.\]

Given the above setup, note that the quantum Grothendieck polynomial $\fg_w^q$ is obtained from the universal one $\fg_w(c)$ simply by specializing the variables $c_p(k)$ to $G_p^k$. It means that (\ref{quiver}) specializes to a formula for $\fg_w^q$ in the obvious way: the left-hand side becomes $\fg_w^q$, and in the right-hand side we replace each variable $c_p(k)$ by $G_p^k$. Recall that $G_p^k$ are computed recursively by Proposition \ref{recg}. This procedure might appear involved, but it is explicit in all of its steps, while no explicit realization of the quantization map is known. Furthermore, we can view (\ref{quiver}) as a reduction formula from the case of an arbitrary permutation $w$ to the case of Grassmannian permutations, i.e. permutations with a unique descent (which correspond to partitions, or Young diagrams).

\section{The quantization map}\label{qmap}

The goal of this section is to give a characterization of the $K$-theory quantization map which is similar to one given in \cite{fgpqsp} for the cohomology quantization map.

We start by recalling from Section \ref{qqa} the quantum Bruhat representation (\ref{actq}) of the quantum quadratic algebra ${\e}_n^q$ on $R\otimes L_n$, the elements $X_k$ defined in (\ref{defxk}), and the fact that they commute, cf. Theorem \ref{comm}.

Consider integers $p,k,n$ with $0\le p\le k< n$. 
We define the polynomial
\begin{equation}\label{deftech}f_p^k:=\frac{1}{1-q_{k-1}}(g_p^{k-1}-g_{p-1}^{k-1})-\frac{q_{k-1}}{1-q_{k-1}}(g_{p-1}^{k-2}-g_{p-2}^{k-2})\,.\end{equation}


\begin{thm}\label{main} Let $w$ be the identity permutation or a permutation in $S_n$ with the first descent in position $k'>k$. We have
\[(1-q_k)(1-X_k)(\fg_w\,f_p^k)=\fg_w\,(g_p^k-g_{p-1}^{k-1})\,.\]
\end{thm}

Note that, by the Pieri formula for Grothendieck polynomials in Theorem \ref{kthm}, the product $\fg_w\,f_p^k$ lies in $L_n$. Also note that, in the classical case (corresponding to $q_i=0$), the theorem reduces to the {\em transition formula} for the polynomials $g_{p}^k$; this formula for Grothendieck polynomials was derived in \cite{lastgp}, cf. also \cite{lenktv}.

The proof of Theorem \ref{main} is postponed to Section \ref{pfmain}. The main idea is to use the Pieri formula for Grothendieck polynomials in Theorem \ref{kthm} in order to expand the product $\fg_w\,f_p^k$ in the basis of Grothendieck polynomials. Indeed, the action of the operator $X_k$ is expressed only in this basis. Therefore, our approach is considerably more complex than the one in \cite{fgpqsp}, in the quantum cohomology case (to be more precise, we refer to the proof of Proposition 5.4 in the mentioned paper). A simpler way to express the action of $X_k$ is not likely to be found because this operator gives rise to paths in the quantum Bruhat graph in which the quantum edges are interspersed with non-quantum edges. Hence it is not possible to separate the ``quantum part'' of $X_k$ from its ``non-quantum part'', as it is done in the quantum cohomology case, cf. the definition (5.2) of the operators ${\mathcal X}_k$ in \cite{fgpqsp}.

\begin{thm}\label{main1} We have $G_p^k(X)(g)=g_p^k\,g$ for any polynomial $g\in L_n^q$ which is symmetric in the variables $x_1,\ldots,x_{k+1}$, $k<n$, and also in the case $g=1$, $k=n$.
\end{thm}

\begin{proof}  
It suffices to consider $g=\fg_w$, where $w$ the identity permutation or a permutation in $S_n$ with the first descent in position $k'>k$. Let us fix such a permutation. We will simultaneously prove the following two relations:
\begin{align}
&G_p^k(X)(\fg_w)=g_p^k\,\fg_w\label{ind1}\\
&(1-q_k)(1-X_k)\left(\gb_p^{k-1}(X)-\gb_{p-1}^{k-1}(X)\right)(\fg_w)=(g_p^k-g_{p-1}^{k-1})\,\fg_w\,.\label{ind2}
\end{align}
We use double induction on $p$ and $k$. The base case $p=0$ is obvious.

Relation (\ref{ind1}) follows easily by combining (\ref{recg3}) with (\ref{ind2}) and the version of (\ref{ind1}) for $p-1$ and $k-1$ (which is part of the induction hypothesis). Indeed, we have
\begin{align*}
G_p^k(X)(\fg_w)&=(1-q_k)(1-X_k)\left(\gb_p^{k-1}(X)-\gb_{p-1}^{k-1}(X)\right)(\fg_w)+\gb_{p-1}^{k-1}(X)(\fg_w)=\\
&=(g_p^k-g_{p-1}^{k-1})\,\fg_w+g_{p-1}^{k-1}\,\fg_w=g_{p}^{k}\,\fg_w\,.\end{align*}

Relation (\ref{ind2}) is derived in several steps. Based on (\ref{recg1}), we first write
\[\gb_p^{k-1}-\gb_{p-1}^{k-1}=G_p^{k-1}+q_{k-1}(1-x_{k-1})(\gb_p^{k-2}-\gb_{p-1}^{k-2})-G_{p-1}^{k-1}-q_{k-1}(1-x_{k-1})(\gb_{p-1}^{k-2}-\gb_{p-2}^{k-2}).\]
Let us now substitute the variables $x_i$ with the operators $X_i$ in the quantum quadratic algebra, and apply both sides of the above equality to $\fg_w$. Based on the induction hypothesis, namely (\ref{ind1}) and (\ref{ind2}) for the pairs $(p,k-1)$ and $(p-1,k-1)$, we obtain
\[\left(\gb_p^{k-1}(X)-\gb_{p-1}^{k-1}(X)\right)(\fg_w)=f_p^k\,\fg_w\,.\]
Relation (\ref{ind2}) for $p$ and $k$ now follows by applying $(1-q_k)(1-X_k)$ to both sides of the last relation and by using Theorem \ref{main}.
\end{proof}

\begin{thm}\label{main2} We have $G_{p_1\ldots p_{n-1}}(X)(1)=g_{p_1\ldots p_{n-1}}$. In particular, $G_p^k(X)(1)=g_p^k$ for any $p\le k\le n$.
\end{thm}

\begin{proof} This is a word-by-word translation of the proof of Theorem 5.5 in \cite{fgpqsp}. In other words, by repeatedly using Proposition \ref{main1}, we have
\begin{align*}
G_{p_1\ldots p_{n-1}}(X)(1)&=G_{p_1\ldots p_{n-2}}(X)(g_{p_{n-1}}^{n-1})=\\
&=G_{p_1\ldots p_{n-3}}(X)(g_{p_{n-2}}^{n-2}g_{p_{n-1}}^{n-1})=\ldots=g_{p_1}^1\ldots g_{p_{n-1}}^{n-1}\,.\end{align*}
Note that, by the Pieri formula in Theorem \ref{kthm}, the expansion of any product $g_{p_i}^i\ldots g_{p_{n-1}}^{n-1}$ in the basis of Grothendieck polynomials contains only Grothendieck polynomials indexed by permutations with first descents in positions greater or equal to $i$. 
\end{proof}

It follows from Theorem \ref{main2} and Proposition \ref{bases} that for every element $f$ of $R\otimes L_n$ there is a unique element $F$ of $R\otimes L_n$ such that $F(X)(1)=f$. Hence, we can define an $R$-linear map $\psi\::\:R\otimes L_n\rightarrow R\otimes L_n$ by $f\mapsto F$. Furthermore, based on the description (\ref{charq}) of the quantization map $\widehat{Q}$ and on the definition of the quantum Grothendieck polynomials $\fg_w^q$, we have the following Corollary. 

\begin{cor}\label{newquant} The map $\psi$ coincides with the quantization map $\widehat{Q}$. In particular, we have 
\[\fg_w^q(X)(1)=\fg_w\,.\]
\end{cor}

\section{Monk-type multiplication formula}\label{monkmul}

 In this section, we derive a Monk-type multiplication formula for our quantum Grothendieck polynomials. This generalizes both the Monk formula for Grothendieck polynomials in \cite{lenktv}, that is, the case $p=1$ of Theorem \ref{kthm}, and the Monk formula for quantum Schubert polynomials in \cite{fgpqsp}, namely Theorem \ref{monkschub}. Before stating the results, let us recall from Section \ref{qqa} the definition of the quantum Bruhat graph of the monomials $q(\pi)$ associated to paths $\pi$ in this graph.

\begin{thm}\label{tmonk} Assume that $w\in S_{n-1}$ and $1\le k<n$. We have
\begin{equation}\label{monk1}x_k\,\fg_w^q=\widehat{Q}(X_k(\fg_w))\,.\end{equation}
In other words, we have (for all $k$):
\begin{equation}\label{monk2}(1-q_k)(1-x_k)\fg_w^q=\sum_\pi (-1)^t {q}(\pi)\fg_{{\rm end}(\pi)}\,;\end{equation}
the summation is over all paths $\pi$  (possibly empty) in the quantum Bruhat graph (of $S_\infty$) of the form
\[
  w=w_0\xrightarrow{\,({a_1,k})\,} 
                    w_{1}\xrightarrow{\,({a_2,k})\,}\ 
           \dotsb\ 
         \xrightarrow{\,({a_s,k})\,}w_{s} 
         \xrightarrow{\,({k,b_1})\,}w_{s+1}
         \xrightarrow{\,({k,b_2})\,}\ \dotsb\ 
         \xrightarrow{\,({k,b_t})\,}w_{s+t}={\rm end}(\pi)
       \,,
\]
where
 \[a_s<a_{s-1}<\dotsb<a_1<k<b_t<b_{t-1}<\dotsb<b_1\,.\]
\end{thm}

\begin{proof}
The conditions on $k$ and $w$ are needed to ensure that $x_k\,\fg_w^q$ lies in $L_n^q$. Let us substitute the variables $x_i$ in (\ref{monk1}) with the operators $X_i$ in the quantum quadratic algebra, and then let both sides act on $1$. By Corollary \ref{newquant}, the right-hand side gives $X_k(\fg_w)$, while the left-hand side gives
\[(X_k\,\fg_w^q(X))(1)=X_k(\fg_w^q(X)(1))=X_k(\fg_w)\,.\]
The explicit formula (\ref{monk2}) now easily follows by recalling the definition (\ref{defxk}) of $X_k$.
\end{proof}

There is a similar formula for expanding the product $\fg_w^q\,\fg_{s_k}^q$. By specializing $q_i=0$, we obtain the Monk formula for Grothendieck polynomials in \cite{lenktv}. Let us first define the {\em quantum $k$-Bruhat graph} as the subgraph of the quantum Bruhat graph whose edges are labeled by pairs $(a,b)$ with $a\le k<b$.

\begin{thm}\label{tmonk2} We have 
\begin{equation}\label{monk3}\fg_w^q\,\fg_{s_k}^q=\sum_\pi (-1)^{\ell(\pi)-1} {q}(\pi)\fg_{{\rm end}(\pi)}\,;\end{equation}
the summation is over all nonempty paths $\pi$ in the quantum $k$-Bruhat graph (of $S_\infty$) of the form
\[
  w=w_0\xrightarrow{\,({a_1,b_1})\,} 
                    w_{1}\xrightarrow{\,({a_2,b_2})\,}\ 
           \dotsb\ 
         \xrightarrow{\,({a_s,b_s})\,}w_{s}={\rm end}(\pi)
       \,,
\]
where 
 \[(a_1,b_1)\prec(a_2,b_2)\prec\ldots\prec(a_s,b_s)\,.\]
This formula has no cancellations.
\end{thm}

\begin{proof}
Let us first recall the explicit form of $\fg_{s_k}^q$ in Corollary \ref{gsk}. We use induction on $k$, where the base case $k=1$ is given by (\ref{monk2}). We have
\begin{equation}\label{indrec}1-\fg_{s_k}^q=\frac{(1-q_k)(1-x_k)}{1-q_{k-1}}(1-\fg_{s_{k-1}}^q)\,.\end{equation}
The multiplication of $\fg_w^q$ by $1-\fg_{s_{k-1}}^q$ is given by the induction hypothesis, where the only differences from (\ref{monk3}) are that the sign in the right-hand side is $(-1)^{\ell(\pi)}$ and the chain $\pi$ might be empty. Let us denote the set of chains in the multiplication formula for $\fg_w^q\,(1-\fg_{s_k}^q)$ by $\Pi_k(w)$, and the set of chains in (\ref{monk2}) by $\Pi_k'(w)$. By combining (\ref{indrec}) with the induction hypothesis and (\ref{monk2}), we can express the product $\fg_w^q\,(1-\fg_{s_k}^q)$ in terms of the following set of concatenated chains:
\[\Pi:=\{\pi_1|\pi_2\::\:\pi_1\in\Pi_{k-1}(w),\;\pi_2\in\Pi_k'({\rm end}(\pi_1))\}\,.\]

Given a generic chain $\pi_1|\pi_2\in \Pi$, let
\[\pi_1=(a_1,b_1),\ldots,(a_s,b_s)\,,\;\;\;\;\pi_2=(k,l_1),\ldots,(k,l_t)\,;\]
here we allow for $l_i<k$; we denote the subchain of $\pi_1$ consisting of transpositions $(i,k)$ by $\pi_1'$ and the subchain of $\pi_2$ consisting of transpositions $(i,k)$ with $i<k$ by $\pi_2'$. Let
\[\overline{\Pi}_1:=\{\pi_1|\pi_2\in\Pi\::\:\pi_1'=\pi_2'=\emptyset\},\;\;\;\overline{\Pi}_2:=\{\pi_1,(k-1,k)|(k-1,k),\pi_2\::\:\pi_1|\pi_2\in \overline{\Pi}_1\}.\]
A chain $\pi_1|\pi_2$ in $\overline{\Pi}_1$ and its obvious pair in $\overline{\Pi}_2$ end in the same permutation, while the contribution of the second chain to the multiplication formula differs from that of the first one by a factor of $-q_{k-1}$. Hence, by summing the two contributions and dividing the result by $1-q_{k-1}$, we obtain $(-1)^{\ell(\pi_1)+\ell(\pi_2)}q(\pi_1|\pi_2)\fg_{{\rm end}(\pi_2)}^q$. Let us also note that, by a simple commutation of transpositions, the chains in $\overline{\Pi}_1$ can be bijected to those in $\Pi_k(w)$.

It remains to prove that the contributions of the chains in $\Pi\setminus(\overline{\Pi}_1\cup \overline{\Pi}_2)$ cancel out. We do this by exhibiting a sign-reversing involution on these chains, such that the contributions cancel out in pairs. This is defined by considering several cases.

{\bf Case 1}: $\pi_1'\ne\emptyset$ (so $b_s=k$) and either $\pi_2'=\emptyset$ or $a_s>l_1$. We map $\pi_1|\pi_2$ to
\[(a_1,b_1),\ldots,(a_{s-1},b_{s-1})|(a_s,b_s),(k,l_1),\ldots,(k,l_t)\,.\]

{\bf Case 2}: $\pi_2'\ne\emptyset$ (so $l_1<k$) and either $\pi_1'=\emptyset$ or $a_s<l_1$. We map $\pi_1|\pi_2$ to
\[(a_1,b_1),\ldots,(a_s,b_s),(l_1,k)|(k,l_2),\ldots,(k,l_t)\,.\]

Chains not yet covered include (in fact, {\em only} include) those for which $\pi_1'\ne\emptyset\ne\pi_2'$ and $(a_s,b_s)=(l_1,k)=(k-1,k)$. Let $\pi_1''$ and $\pi_2''$ be the subchains of $\pi_1'$ and $\pi_2'$ obtained by removing the transposition $(k-1,k)$.

{\bf Case 3}: $\pi_1''\ne\emptyset$ (so $b_{s-1}=k$) and either $\pi_2''=\emptyset$ or $a_{s-1}>l_2$. We map $\pi_1|\pi_2$ to
\[(a_1,b_1),\ldots,(a_{s-2},b_{s-2}),(k-1,k)|(k-1,k),(a_{s-1},b_{s-1}),(k,l_2),\ldots,(k,l_t)\,.\]

{\bf Case 4}: $\pi_2''\ne\emptyset$ (so $l_2<k$) and either $\pi_1''=\emptyset$ or $a_{s-1}<l_2$. We map $\pi_1|\pi_2$ to
\[(a_1,b_1),\ldots,(a_{s-1},b_{s-1}),(l_2,k),(k-1,k)|(k-1,k),(k,l_3),\ldots,(k,l_t)\,.\]

It is easy to check that the definitions are correct, and the corresponding contributions cancel out. Furthermore, the above cases exhaust all possibilities because, as it is easy to check, we can have neither $(a_s,b_s)=(l_1,k)\ne(k-1,k)$, nor simultaneously $(a_s,b_s)=(l_1,k)=(k-1,k)$ and $(a_{s-1},b_{s-1})=(l_2,k)$.
\end{proof}

Based on Theorem \ref{tmonk2}, we conjecture a Pieri-type formula for quantum Grothendieck polynomials. By analogy with Definition \ref{Pieri_chain} for $k$-Pieri chains, we first  define a {\em quantum $k$-Pieri chain} as a path $\pi$ of the form (\ref{chain}) in the quantum $k$-Bruhat graph which satisfies the same conditions (P1) and (P2). A marking of a quantum $k$-Pieri chain is defined similarly (by conditions (M1)-(M3)); the same is true for the coefficient $m_p(\pi)$. In addition, we need to keep track of the down-steps in $k$-Bruhat order, so we define $m^q_p(\pi):=m_p(\pi)q(\pi)$, where $q(\pi)$ was defined in Section \ref{qqa}. With this notation, we can state our conjecture, which is a common generalization of the Pieri formulas for Grothendieck and quantum Schubert polynomials in \cite{lasptf} and \cite{posqvp}, respectively.

\begin{conj}\label{qkpieri}
  We have that
 \[
   \fg_w^q\,G_p^k=
   \sum_\pi m_p^q(\gamma)\fg_{{\rm end}(\pi)}^q\,,
 \]
 where the sum is over all quantum $k$-Pieri chains $\pi$ (on the infinite symmetric group) 
 that begin at $w$.
 This formula has no cancellations.
\end{conj}

\section{Main conjecture and applications}\label{mainconj}

We now state our main conjecture. This is the analog of Theorem \ref{repschub}, which states that quantum Schubert polynomials represent Schubert classes in the quantum cohomology of $Fl_n$.

\begin{conj}\label{maincj} The quantum Grothendieck polynomials $\fg_w^q$ are representatives for Schubert classes in $QK(Fl_n)\simeq\ZZ[q,x]/\widehat{I}_n^q$. \end{conj}

The proof of Theorem \ref{repschub} in \cite{fgpqsp} has a single geometric component, which is a result in \cite{fonqcf}, stating that the polynomial $E_p^k$ represents the Schubert class indexed by the permutation $c[k,p]=(k-p+1,k-p+2,\ldots,k+1)$. A proof of the above conjecture would most probably still require one to prove geometrically that the polynomial $\widehat{E}_p^k$ represents the Schubert class in $QK(Fl_n)$ indexed by $c[k,p]$. However, as explained at the end of this section, more geometric information is needed. 

There is strong evidence for the above conjecture. First of all, the polynomials $\eb_p^k$ defined in (\ref{defepk}) provide a link between $QK(Fl_n)$ and quantum Grothendieck polynomials. Indeed, a presentation for the former (Theorem \ref{presqk}) and the definition of the latter (Definition \ref{defqgp}) are both given in terms of the polynomials $\eb_p^k$. Secondly, the operators $X_k$ in the quantum quadratic algebra $\e_n^q$ provide another link between $QK(Fl_n)$ and quantum Grothendieck polynomials. Indeed, on the one hand, the subalgebra generated by these operators inside $\e_n^q$ is isomorphic to $QK(Fl_n)$, as stated in Theorem \ref{qkenq}; on the other hand, the operators $X_k$ realize the multiplication of $\fg_w^q$ by the variable $x_k$, as stated in Theorem \ref{tmonk}. Thirdly, our Monk formula for quantum Grothendieck polynomials in Theorem \ref{tmonk2} is the natural common generalization of the similar formulas for Grothendieck polynomials \cite{lenktv} (i.e., the case $p=1$ of Theorem \ref{kthm}) and quantum Schubert polynomials \cite{fgpqsp} (i.e., Theorem \ref{monkschub}); both of the latter formulas are multiplication formulas for the corresponding Schubert classes. Furthermore, our experiments indicate that the polynomials $\fg_w^q$ in this paper are the unique family of polynomials satisfying the Monk-type formula in Theorem \ref{tmonk2}.

Combined with Theorem \ref{tmonk2}, the above conjecture would confirm the type $A$ version of the Monk-type multiplication formula for Schubert classes in quantum $K$-theory that was conjectured in \cite[Section 17]{lapawg}.

Let us now recall the definition of the product in the ring $QK(Fl_n)$. 
Given a collection of nonnegative integers $d=(d_1,\dots,d_r)$,
called multidegree, we define $q^d$ as in (\ref{notexp}).  
As a $\ZZ[q]$-module, the quantum $K$-theory is defined as
$QK(Fl_n):= K(Fl_n)\otimes_\ZZ \ZZ[q]$. 
Let $[w]$ denote the class of the structure sheaf of the Schubert variety
$X_{w}$.
Then the classes of $[w]$ form a $\ZZ[q]$-basis of $QK(Fl_n)$. 
The multiplication in $QK(Fl_n)$ is a deformation of the classical
multiplication:
\begin{equation}\label{qkmult}
[u]\circ [v]=\sum_d q^d\sum_{w\in S_n} N_{uv}^w(d) \, [w]\,,
\end{equation}
where the first sum is over all multidegrees $d$, and $N_{uv}^w(d)$ is the 
{\it
$3$-point quantum $K$-invariant of Gromov-Witten type} for $[u]$, $[v]$, and the
quantum dual of $[w]$ (see the discussion at the end of this section for more on duality in quantum $K$-theory). As defined in \cite{leeqkt}, this invariant is
the $K$-theoretic push-forward to $\mathrm{Spec}\, {\mathbb C}$ of some natural vector
bundle on the moduli space $\overline{M}_{3,0}(Fl_n,d)$ (via the orientation
defined by the virtual structure sheaf). 
The associativity of the quantum $K$-product
was established in \cite{leeqkt}, based on a sheaf-theoretic version of
an argument of WDVV-type. 

Assuming that Conjecture \ref{maincj} is true, we can compute the quantum $K$-invariants $N_{uv}^w(d)$ by expanding the product $F:=\fg_u^q\fg_v^q$ in  the basis of quantum Grothendieck polynomials (cf. Theorem \ref{tmonk2} and Conjecture \ref{qkpieri}). This is realized by Algorithm \ref{expgro}. 

\begin{rem} (1) The algorithm terminates, that is, the product $\fg_u^q\fg_v^q$ has an expansion in the basis of quantum Grothendieck polynomials with coefficients in $\ZZ[q]$, due to the geometric reasons mentioned above. However, from a purely algebraic point of view, this is far from clear beyond the Monk-type formula in Theorem \ref{tmonk2}. 

(2) The algorithm can be stopped at any time, and the conjectured computation of the quantum $K$-invariants obtained so far is finished. Indeed, this computation is not continued by subsequent iterations in the algorithm due to the increasing condition on the list $L$ that the algorithm outputs. 
\end{rem}

\begin{exa}\label{qkinvgw}{\rm Let us compute the expansion of $\fg_{321}^q\fg_{231}^q$, to which the Monk-type formula in Theorem \ref{tmonk2} does not apply. The two quantum Grothendieck polynomials are found in Example \ref{exs3}. We will also need some quantum Grothendieck polynomials for $S_4$, which were found by a computer. Steps 2 and 3 in the first iteration of the algorithm provide
\[F_0=x_1^3x_2^2+q_1^2x_1+2q_1x_1^2x_2=\fs^q_{4312}+q_2\fs^q_{4123}+q_1q_2\fs^q_{132}\,.\]
Steps 2 and 3 in the second iteration provide
\[F_0=q_2(-x_1^3x_2-x_1^3x_3-2q_1x_1^2-q_1x_1x_2+q_1x_1x_3-q_1x_2^2+q_1^2+q_1q_2)=-q_2\fs^q_{4132}-q_1q_2\fs^q_{1342}-q_1q_2\fs^q_{1423}.\]
Steps 2 and 3 in the third iteration provide
\[F_0=q_1q_2(x_1^2x_2+x_1^2x_3+x_1x_2^2+x_1x_2x_3+x_2^2x_3+q_1x_1+q_1x_2-q_1x_3+q_2x_2)=q_1q_2\fs^q_{1432}\,.\]
The algorithm stops after the third iteration. Hence we have
\[\fg_{321}^q\fg_{231}^q=(\fg^q_{4312}+q_2\fg^q_{4123}+q_1q_2\fg^q_{132})-(q_2\fg^q_{4132}+q_1q_2\fg^q_{1342}+q_1q_2\fg^q_{1423})+q_1q_2\fg^q_{1432}\,.\]
The first bracket gives the Gromov-Witten invariants in quantum cohomology (e.g., see \cite[Section 2.3]{fgpqsp}); indeed, the expansion of $F_0$ in the first iteration is precisely the expansion of $\fs_{321}^q\fs_{231}^q$. Classical $K$-theory gives $\fg_{321}\fg_{231}=\fg_{4312}$. So, starting with the second bracket, we have information which, conjecturally, is given only by quantum $K$-theory. For instance, conjecturally, we have $N_{3214,2314}^{1432}(1,1)=1$ for $Fl_4$. 
}\end{exa}

Brion \cite{bripgg} proved that the structure constants of the $K$-theory of a generalized flag variety $G/B$ have alternating signs. Based on this results, as well as on our Monk-type formula (Theorem \ref{tmonk2}), we make the following conjecture, which is also supported by Example \ref{qkinvgw}.

\begin{conj}\label{signalt} The quantum $K$-invariants of Gromov-Witten type for $Fl_n$ have alternating signs, i.e., we have
\[(-1)^{\sum_i d_i+\ell(w)-\ell(u)-\ell(v)}N_{uv}^w(d)\ge 0\,.\]
The same result holds for a generalized flag variety $G/B$. 
\end{conj}

We conclude this section with a discussion of the pairing in quantum $K$-theory \cite{leeqkt}, which is a deformation of the natural pairing $\chi$ in $K$-theory (given by Euler characteristic of vector bundles). More precisely, one defines the quantum $K$-theory pairing on the Schubert classes $[u]$ and $[v]$ by
\begin{equation}\label{pair}\langle\!\langle[u],[v]\rangle\!\rangle:=\chi([u][v])+\sum_d q^d N_{uv}(d)\,;\end{equation}
here $N_{uv}(d)=N_{uvw}(d)$ for $w={\rm id}$ are the 2-point quantum $K$-invariants for $[u]$ and $[v]$. Due to this deformation, we are not able to define a pairing on quantum Grothendieck polynomials purely algebraically. Indeed, more geometric information is needed in order to compute the 2-point quantum $K$-invariants in (\ref{pair}). In fact, the 3-point quantum $K$-invariants $N_{uvw}(d)$ of $[u]$, $[v]$, $[w]$, and the quantum $K$-invariants $N_{uv}^w(d)$ in (\ref{qkmult}) determine each other, but only  in the presence of a {\em metric}, i.e., of the 2-point invariants $N_{uv}(d)$. Recall that it is the $N_{uv}^w(d)$ which, conjecturally, can be computed purely based on quantum Grothendieck polynomials. 

The above situation is drastically different from quantum cohomology. Indeed, the pairing $(\!(f,g)\!)$ on the algebra $\ZZ[q,x]/I_n^q$ (which is isomorphic to $QH^*(Fl_n)$ by Theorem \ref{qcfln}) is defined simply as the coefficient of $\fs_{w_\circ}^q$ in the expansion of $fg$ in the basis of quantum Schubert polynomials. (Alternatively, $(\!(f,g)\!)$ is the coefficient of the staircase monomial $x_1^{n-1}x_2^{n-2}\ldots x_{n-1}$ in the monomial expansion of $fg$ in $\ZZ[q,x]/I_n^q$.) This simpler situation allowed Fomin, Gelfand, and Postnikov to prove the orthogonality of quantum Schubert polynomials purely algebraically in \cite{fgpqsp}. Then, based on this information and a small amount of geometric information in \cite{fonqcf}, they proved that the quantum Schubert polynomials represent Schubert classes in $QH^*(Fl_n)$. Unfortunately, as explained above, this approach does not work for quantum Grothendieck polynomials. Hence, the larger complexity of quantum $K$-theory requires more geometric information in order to prove Conjecture \ref{maincj}.

\section{Quantum double Grothendieck polynomials}\label{qdgroth}

We start by recalling the {\em double Grothendieck polynomials} $\fg_w(x,y),$ $w\in S_n,$ which were defined by Lascoux and Sch\"utzenberger \cite{lasssa}; they represent Schubert classes in the equivariant $K$-theory of $Fl_n$. Let 
\[ \fg_{w_\circ}(x,y)=\prod_{i+j\leq n}(x_i+y_j-x_i y_j) \] 
be the double Grothendieck polynomial for the longest element $w_\circ \in S_n.$ 
The double Grothendieck polynomial for an element $w\in S_n$ 
is obtained by applying the isobaric divided difference operators $\pi_i,$ 
where $1\le i\le n-1$, to the polynomial $\fg_{w_\circ}(x,y)$. 
The isobaric divided difference operator $\pi_i$ is by definition 
\begin{equation}\label{isodivdiff} \pi_i =1+(1-x_i) \partial_i\,, \end{equation}
where $\partial_i$ is the divided difference operator in (\ref{divdiff}).  
If $w\in S_n$ has a reduced decomposition $w=s_{i_1}\ldots s_{i_l},$ 
the operator $\pi_w$ is defined by $\pi_w:=\pi_{i_1}\ldots \pi_{i_l}$;  
this definition is correct, i.e., is independent of the choice of the reduced decomposition, 
because the operators $\pi_i$ satisfy the braid relations. 
The double Grothendieck polynomial $\fg_w(x,y)$ is defined 
by the formula 
\[ \fg_w (x,y):= \pi_{w^{-1}w_\circ}^{(x)}\fg_{w_\circ}(x,y). \] 
Here the isobaric divided difference operator $\pi_{w^{-1}w_\circ}^{(x)}$ 
acts on the $x$ variables only. 

In order to define quantum double Grothendieck polynomials, note first that 
\[ \fg_{w_\circ}(x,y)=\prod_{i=1}^{n-1}\left( 
1+\sum_{j=1}^i (-1)^j(1-x_{n-i})^j f_j^i(y) \right) \,, \] 
where $f_j^i(y):=e_j(1-y_1,\ldots,1-y_i)$, as defined in (\ref{defsfpk}). 

\begin{dfn}\label{qdgro} 
The polynomial $\fg^q_{w_\circ}(x,y)$ for the element $w_\circ \in S_n$ 
is given by the formula 
\[ \fg^q_{w_\circ}(x,y):=\prod_{i=1}^{n-1}\left( 
1+\sum_{j=1}^i (-1)^j(1-x_{n-i})^j F^i_j(y) \right) \,, \] 
where the polynomial $F^i_j$ is defined in {\rm (\ref{deffpk})}. 
For an arbitrary element $w\in S_n,$ the {\em quantum 
double Grothendieck polynomial} $\fg^q_w(x,y)$ is 
\[ \fg^q_w(x,y):= \pi_{w^{-1}w_\circ}^{(x)}\fg^q_{w_\circ}(x,y)\,. \] 
\end{dfn}

\begin{rema}
We clearly have $\fg_w^q(x,y)|_{q=0}=\fg_w(x,y)$.
\end{rema}

Lascoux and Sch\"utzenberger~\cite{lassfm} defined the {\em dual Grothendieck polynomials} $\fh_w$, for $w\in S_n$, by 
\begin{equation}\label{hpoly} \fh_w=\sum_{v\in S_n,\,v \geq w} (-1)^{\ell(v)-\ell(w)}\fg_v\, . \end{equation}
Unlike Grothendieck polynomials, these polynomials are unstable, so they depend on $n$. They represent the $K$-theory classes dual to the classes of structure sheaves with respect to the natural intersection pairing in $K$-theory, see \cite[Proposition 2.1]{lrsgpp}. In fact, it was shown by Brion and Lakshmibai \cite{balgas} that a dual class corresponds to the ideal sheaf of the boundary $X_w-X_w^{\circ}$ of the Schubert variety $X_w$. Several combinatorial formulas for the dual Grothendieck polynomials can be found in \cite[Section 6]{lrsgpp}. 

Let us now recall the Cauchy 
identity for the classical Grothendieck polynomials, which is due to Fomin and Kirillov \cite{fakgpy} 
(see also \cite[Proposition 2]{kirqgp}): 
\begin{equation}\label{grocauchyid} \fg_{w_\circ}(x,y)= \sum_{w \in S_n} \fh_w \fg_{ww_\circ}(y) \,. \end{equation}This identity is generalized as follows for the quantum Grothendieck 
polynomials. 

\begin{thm}\label{cauchy} We have
\begin{equation}\label{cauchyid} \fg^q_{w_\circ}(x,y)= \sum_{w \in S_n} \fh_w \fg^q_{ww_\circ}(y)\,. \end{equation}
\end{thm}

\begin{proof}
Let us consider the $K$-theoretic quantization map $\widehat{Q}^{(y)}$ 
with respect to the $y$ variables as a $\ZZ[q,x]$-linear map. 
By applying this map to both sides of the Cauchy 
identity for the classical Grothendieck polynomials (\ref{grocauchyid}), we have
\[ \widehat{Q}^{(y)}(\fg_{w_\circ}(x,y))= 
\widehat{Q}^{(y)}\left(\sum_{w \in S_n} \fh_w \fg_{ww_\circ}(y)\right)= 
\sum_{w \in S_n} \fh_w \fg^q_{ww_\circ}(y) \,. \] 
The proof is concluded by the following calculation:
\begin{align*}
\widehat{Q}^{(y)}(\fg_{w_\circ}(x,y))&=\widehat{Q}^{(y)}\left( 
\prod_{i=1}^{n-1}\left( 
1+\sum_{j=1}^i (-1)^j(1-x_{n-i})^j f_j^i(y) 
\right) \right) \\ 
 &= \prod_{i=1}^{n-1}\left( 
1+\sum_{j=1}^i (-1)^j(1-x_{n-i})^j \widehat{Q}^{(y)}(f_j^i(y)) 
\right) \\ 
 &= \prod_{i=1}^{n-1}\left( 
1+\sum_{j=1}^i (-1)^j(1-x_{n-i})^j F^i_j(y) \right)=\fg^q_{w_\circ}(x,y)\,. 
\end{align*}
The second and third equalities follow from Proposition \ref{quantf}.
\end{proof}

We now derive a corollary of the Cauchy identity, which leads to an explicit recursive construction of the quantum Grothendieck polynomials. We need the following lemma.

\begin{lem}\label{actp}
We have $\pi_w(\fg_v)=1$ if and only if $v\le w$ in the Bruhat order on $S_n$.
\end{lem}

\begin{proof}
It is well-known that, if $\ell(ws_i)>\ell(w)$, where $s_i$ is the adjacent transposition $t_{i,i+1}$, then
\begin{equation}\label{isdd} \pi_{ws_i}(\fg_v)=\case{\pi_w(\fg_{vs_i})}{\ell(vs_i)<\ell(v)}{\pi_w(\fg_v)}\end{equation}
We proceed by induction on $\ell(w)$, given a fixed permutation $v$. Clearly, we can have $\pi_w(\fg_v)=1$ only if $\ell(w)\ge\ell(v)$. So induction starts at $\ell(w)=\ell(v)$, in which case the statement is easily checked. Now assume that the statement is known for all permutation $w$ of a fixed length $k$. Pick a permutation of length $k+1$, which can be written as $ws_i$, where $\ell(w)=k$; in particular, we have $\ell(ws_i)>\ell(w)$.  
The induction step is completed based on (\ref{isdd}), the induction hypothesis, and the following recursive characterization of the Bruhat order on any Coxeter group (e.g., see \cite[Theorem 1.1]{deoscb} or \cite{vermib}): if $\ell(ws_i)>\ell(w)$, then
\[v\le ws_i\;\;\Longleftrightarrow\;\;\casefour{\ell(vs_i)<\ell(v)}{vs_i\le w}{\ell(vs_i)>\ell(v)}{v\le w}\]
\end{proof}

\begin{cor}\label{d2s}
We have 
\[  \fg_w^q=\fg_{w^{-1}}^q(y,x)|_{y=0}\,.\] 
\end{cor}

\begin{proof}
Let us apply the operator $\pi_{w^{-1}w_\circ}^{(x)}$ to both sides of (\ref{cauchyid}). By Definition \ref{qdgro}, we obtain $\fg_w^q(x,y)$ on the left-hand side. By (\ref{hpoly}), the right-hand side of (\ref{cauchyid}) can be rewritten as
\[ \sum_{u \in S_n} \fg^q_{uw_\circ}(y) \left(\sum_{v \geq u} (-1)^{\ell(v)-\ell(u)}\fg_v\right)\,.\]
By applying the operator $\pi_{w^{-1}w_\circ}^{(x)}$ to this expression and then setting the $x$ variables to 0, we obtain the following expression, based on Lemma \ref{actp}:
\begin{equation}\label{mobbru} \sum_{u \le w^{-1}w_\circ} \fg^q_{uw_\circ}(y) \left(\sum_{u\le v\le w^{-1}w_\circ} (-1)^{\ell(v)-\ell(u)}\right)\,.\end{equation}
Now let us recall the M\"obius function of the Bruhat order on a Coxeter group, which is the unique integer function $\mu$ on pairs $u\le v$ in the group such that $\mu(v,v)=1$ and $\sum_{u\le x\le v}\mu(u,x)=0$ if $u<v$. It is a classical result of Verma \cite{vermib} that $\mu(u,v)=(-1)^{\ell(v)-\ell(u)}$ for all $u\le v$. This implies that the interior summation in (\ref{mobbru}) is 0 unless $u=w^{-1}w_\circ$. Hence, the expression in (\ref{mobbru}) is simply $\fg_{w^{-1}}^q(y)$. 
\end{proof}

\begin{rema}
Clearly, Theorem \ref{cauchy} generalizes the Cauchy identity for Grothendieck polynomials (\ref{grocauchyid}), as well as the Cauchy identity for quantum Schubert polynomials in Theorem \ref{schucauchy}. As far as Corollary \ref{d2s} is concerned, upon setting the $q$ variables to 0 in it, we obtain the well-known relationship between double and ordinary Grothendieck polynomials $\fg_w=\fg_w(x,y)|_{y=0}$. Indeed, it is known that $\fg_{w^{-1}}(y,x)=\fg_w(x,y)$. Furthermore, Corollary \ref{d2s} also extends Theorem \ref{schud2s} related to quantum Schubert polynomials.
\end{rema}

\section{Proof of Theorem \ref{main}}\label{pfmain}

Recall the definition (\ref{deftech}) of the polynomials $f_p^k$, for $0\le p\le k< n$. In order to find the expansion of $\fg_w\,f_p^k$ in the basis of Grothendieck polynomials, we use the Pieri formula for Grothendieck polynomials in Theorem \ref{kthm}. 

For simplicity, we consider first the case when $w$ is a Grassmannian permutation; in other words, $w$ has a unique descent in position $k'>k$. We need the special case of Theorem \ref{kthm} corresponding to the products $\fg_w\,g_r^l$ for $r\in\{p-2,p-1,p\}$ and $l\in\{k-2,k-1,k\}$, where $p,k$ are the ones in Theorem \ref{main}. 

We denote by $\Gamma_l$ the set of $l$-Pieri chains starting at $w$. Such a chain $\gamma$ has the form
\begin{align}\label{pierichain}&(a_1,b_1),(a_1-1,b_1),\ldots,(a_1',b_1),\ldots, (a_m,b_m),(a_m-1,b_m),\ldots,(a_m',b_m),\\&(l,l+1),(l-1,l+1),\ldots,(l+1-h,l+1)\,,\nonumber\end{align}
where $h,m\ge 0$ and
\begin{equation}\label{ab}n\ge b_1>\ldots>b_m>k'\,,\;\;\;\min(l,l+1-h)\ge a_1\ge a_1'>\ldots>a_m\ge a_m'\ge 1\,.\end{equation}
Let us note that 
\begin{equation}\label{ab1}w(j+1)=w(j)+1\;\;\;\mbox{for}\;\;\;a_i'\le j<a_i\,,\;\;\;\mbox{and}\;\;\;w(b_i)=w(a_i)+1\,,\end{equation}
for $i=1,\ldots,m$. Let $A:=\sum_{i=1}^m (a_i-a_i')$ if $m\ge 1$, and $A:=0$ if $m=0$. We will use the notion of an $i$-subchain of $\gamma$ (for a fixed $i$), which is the subchain consisting of all covers labeled $(\,\cdot\,,i)$. The first such subchain will be called the initial subchain. 

It is not hard to see that we have
\begin{equation}\label{coeff}
m_r(\gamma)=\casethree{(-1)^{m+A+h-r}\binom{m}{r-A-h}}{l+1-h>a_1\;\;\mbox{and}\;\;h\ge 1}{(-1)^{m+A-r}\binom{m-1}{r-A-1}}{h=0}{(-1)^{m+A+h-r}\binom{m-1}{r-A-h}}{l+1-h=a_1}
\end{equation}
In the last case, we obviously have $h,m\ge 1$. In order to verify (\ref{coeff}), note first that the marking rules (M1)-(M3) require us to mark all covers labeled $(l,l+1),\ldots,(l+2-h,l+1)$, as well as those labeled $(\,\cdot\,,b_i)$ for $i=1,\ldots,m$, with the exception of $(a_2',b_2),\ldots,(a_m',b_m)$. Furthermore, if $m=0$ we must also mark the cover labeled $(l+1-h,l+1)$, while if $l+1-h=a_1$ then the mentioned cover must not be marked. All other covers may or may not be marked, and hence (\ref{coeff}) follows. This formula will be used both explicitly and implicitly several times below.

Given $l\in\{k-2,k-1\}$, let us denote by $\widehat{\Gamma}_l$ the set of all concatenations $\gamma=\pi|\mu$ of a chain $\pi$ in $\Gamma_l$ and a chain $\mu$ in the quantum Bruhat graph, possibly empty, of the following form:
\begin{align*}
  &{\rm end}(\pi)=w_0\xrightarrow{\,({c_1,k})\,} 
                    w_{1}\xrightarrow{\,({c_2,k})\,}\ 
           \dotsb\ 
         \xrightarrow{\,({c_s,k})\,}w_{s} 
         \xrightarrow{\,({k,d_1})\,}w_{s+1}
         \xrightarrow{\,({k,d_2})\,}\ \dotsb\ 
         \xrightarrow{\,({k,d_t})\,}w_{s+t} 
       \,,\\[0.05in]
    & \mbox{where }\;\;1\le c_s<c_{s-1}<\dotsb<c_1<k<d_t<d_{t-1}<\dotsb<d_1\le n\,.\end{align*}
We call $\mu$ a {\em Monk} (sub){\em chain} (of $\gamma$) since our Monk formula for the quantum Grothendieck polynomials is expressed in terms of such chains. We also define the weight of an edge labeled $(i,j)$ to be $1$ if it corresponds to an increase in length by $1$, and $q_i\ldots q_{j-1}$ otherwise. We then define $q(\mu)$ to be $(-1)^t$ times the product of the weights of all edges, and $m_r(\gamma):=m_r(\pi)q(\mu)$. 

Based on the Pieri formula and on the action of the operator $X_k$, the identity in Theorem \ref{main} is equivalent to
\begin{align}\label{toprove}&\sum_{\gamma\in\widehat{\Gamma}_{k-1}} \frac{1}{1-q_{k-1}}m_p(\gamma)\fg_{\eg} + \sum_{\gamma\in\widehat{\Gamma}_{k-1}} -\frac{1}{1-q_{k-1}}m_{p-1}(\gamma)\fg_{\eg}+ \\+&\sum_{\gamma\in\widehat{\Gamma}_{k-2}} -\frac{q_{k-1}}{1-q_{k-1}}m_{p-1}(\gamma)\fg_{\eg} +  \sum_{\gamma\in\widehat{\Gamma}_{k-2}} \frac{q_{k-1}}{1-q_{k-1}}m_{p-2}(\gamma)\fg_{\eg} =\nonumber \\=&\sum_{\gamma\in\Gamma_k}m_{p}(\gamma)\fg_{\eg}+\sum_{\gamma\in\Gamma_{k-1}}-m_{p-1}(\gamma)\fg_{\eg}\,.\nonumber\end{align}
We will partition the nonzero terms in the four sums on the left-hand side of (\ref{toprove}) into blocks such that: (i) the sum of the terms in some blocks is 0; (ii) each of the remaining blocks is paired up with one or two terms on the right-hand side of (\ref{toprove}) such that the corresponding sums are identical. 

We now describe the mentioned blocks in terms of the chains corresponding to them. Since there are several types of such blocks, we will consider several cases. Certain chains will be used several times below, so we introduce them now. Let $\gamma_i$ be the chain in $\widehat{\Gamma}_{k-1}$ of the following form (cf. (\ref{pierichain})):
\[(a_1,b_1),\ldots,(a_m,b_m'),(k-1,k),\ldots,(k-i+1,k)|(k-i,k),\ldots,(k-h,k)\,,\]
where $m\ge 0$, $h\ge 1$, $i=1,\ldots,h+1$, and (\ref{ab}) holds with $l=k-1$. Let $\delta_i$ be the chain in $\widehat{\Gamma}_{k-2}$ of the form
\[(a_1,b_1),\ldots,(a_m,b_m'),(k-2,k-1),\ldots,(k-i,k-1)|(k-1,k),(k-i-1,k),\ldots,(k-h,k)\,,\]
where $m\ge 0$, $h\ge 1$, $i=1,\ldots,h$, (\ref{ab}) holds with $l=k-1$, and $a_1\le k-2$. It is easy to see that all the chains $\gamma_i$ and $\delta_i$ end in the same permutation (just commute the transposition $(k-1,k)$ in ${\delta}_i$ past $i-1$ transpositions to its left). Let $p':=p-A-1$, where $A$ is defined as above; this notation will be used throughout this section.

{\bf Case 1.} This case corresponds to blocks of terms on the left-hand side of (\ref{toprove}) which are matched with terms on the right-hand side. We will use the chains $\gamma_i$ and $\delta_i$ introduced above. The crucial assumption we make throughout this case, without mentioning it again, is that $h\ge 1$ and $a_1<k-h$.

{\bf Case 1.1.} This case corresponds to the chains $\gamma_i$ in $\widehat{\Gamma}_{k-1}$. Based on (\ref{coeff}), we have 
\begin{align*}\sum_{i=1}^{h+1}m_p(\gamma_i)&=(-1)^{m-p'-1}\binom{m-1}{p'}+(-1)^{m-p'}\binom{m}{p'}+(-1)^{m+1-p'}\binom{m}{p'-1}+\\&+\ldots+(-1)^{m+h-p'-1}\binom{m}{p'-h+1}\,,\\
\sum_{i=1}^{h+1}m_{p-1}(\gamma_i)&=(-1)^{m-p'}\binom{m-1}{p'-1}+(-1)^{m+1-p'}\binom{m}{p'-1}+\ldots+\\&+(-1)^{m+h-p'-1}\binom{m}{p'-h+1}+(-1)^{m+h-p'}\binom{m}{p'-h}\,.
\end{align*}
Hence, by using Pascal's identity, we have 
\[\sum_{i=1}^{h+1}m_p(\gamma_i)-\sum_{i=1}^{h+1}m_{p-1}(\gamma_i)=(-1)^{m+h-p'-1}\binom{m}{p'-h}\,.\]

This case also corresponds to the chains $\delta_i$ in $\widehat{\Gamma}_{k-2}$. By a similar calculation to the one above, we have
\[-\sum_{i=1}^{h}m_{p-1}(\delta_i)+\sum_{i=1}^{h}m_{p-2}(\delta_i)=(-1)^{m+h-p'}\binom{m}{p'-h}\,.\]
We conclude that the sum of the coefficients in the terms on the left-hand side of (\ref{toprove}) corresponding to the chains considered above is $(-1)^{m+h-p'-1}\binom{m}{p'-h}$. The chains $\gamma_i$ can be viewed as a single chain $\gamma$ in $\Gamma_{k-1}$. The coefficient in the corresponding term on the right-hand side of (\ref{toprove}), namely $-m_{p-1}(\gamma)$, is also $(-1)^{m+h-p'-1}\binom{m}{p'-h}$.

{\bf Case 1.2.} This case is very similar to the previous one. We consider the chains 
\[{\gamma}'_i:=\gamma_i,(k,k+1)\]
 in $\widehat{\Gamma}_{k-1}$, for $i=1,\ldots,h+1$, and 
\[{\delta}_i':=\delta_i,(k,k+1)\]
 in $\widehat{\Gamma}_{k-2}$, for $i=1,\ldots,h$. Based on the above calculation, the sum of the coefficients in the terms on the left-hand side of (\ref{toprove}) corresponding to the chains $\gamma_i'$ and $\delta_i'$ is $(-1)^{m+h-p'}\binom{m}{p'-h}$. Note that the sign is different from the one in Case 1.1 because each Monk subchain $\mu$ now has one step of the form $(k,\,\cdot\,)$, so $q(\mu)=-1$. The corresponding term on the right-hand side of (\ref{toprove}) is given by the chain 
\[\gamma'=(a_1,b_1),\ldots,(a_m,b_m'),(k,k+1),(k-1,k+1),\ldots,(k-h,k+1)\]
in $\Gamma_k$. Note that the ends of the chains $\gamma_i'$, $\delta_i'$, and $\gamma'$ coincide (indeed, commute the transposition $(k,k+1)$ in $\gamma'$ past the transpositions to its right). Furthermore, we have $m_p(\gamma')=(-1)^{m+h-p'}\binom{m}{p'-h}$. 

{\bf Case 1.3.} This case is again similar to Case 1.1. We consider the chains 
\[{\gamma}''_i:=\gamma_i,(k,b_0)\]
 in $\widehat{\Gamma}_{k-1}$, for $i=1,\ldots,h+1$, and 
\[{\delta}_i'':=\delta_i,(k,b_0)\]
 in $\widehat{\Gamma}_{k-2}$, for $i=1,\ldots,h$, where $b_0>k'$, and, in fact, $b_0>b_1$ if $m\ge 1$. We also assume that there is $h'$ with $1\le h'\le h$, such that
\begin{equation}\label{case13}w(j+1)=w(j)+1\;\;\mbox{for}\;\;k-h\le j<k-h'\,,\;\;\mbox{and}\;\;w(b_0)=w(k-h')+1\,.\end{equation}
Like in Case 1.2, the sum of the coefficients in the terms on the left-hand side of (\ref{toprove}) corresponding to the chains $\gamma_i''$ and $\delta_i''$ is $(-1)^{m+h-p'}\binom{m}{p'-h}$. The corresponding term on the right-hand side of (\ref{toprove}) is given by the chain 
\[\gamma''=(k-h',b_0),\ldots,(k-h,b_0),(a_1,b_1),\ldots,(a_m,b_m'),(k-1,k),\ldots,(k-h',k)\]
in $\Gamma_{k-1}$. This chain ends in the same permutation as $\gamma_i''$ and $\delta_i''$; indeed, we can easily compute the composition of the transpositions below:
\begin{align*}(k-1,k)\ldots(k-h,k)(k,b_0)&=(k-h,\ldots,k-1,k,b_0)=\\
&=(k-h',b_0)\ldots(k-h,b_0)(k-1,k)\ldots(k-h',k)\,.\end{align*}
Also note that $-m_{p-1}(\gamma'')=(-1)^{m+h-p'}\binom{m}{p'-h}$; here, as opposed to the cases above, the last case in (\ref{coeff}) was used.

{\bf Case 1.4.} This case is a combination of Cases 1.2 and 1.3. We consider the chains
\[{\gamma}'''_i:=\gamma_i,(k,b_0),(k,k+1)\]
 in $\widehat{\Gamma}_{k-1}$, for $i=1,\ldots,h+1$, and 
\[{\delta}_i''':=\delta_i,(k,b_0),(k,k+1)\]
 in $\widehat{\Gamma}_{k-2}$, for $i=1,\ldots,h$, under the same conditions as above. The corresponding term on the right-hand side of (\ref{toprove}) is given by the chain 
\[(k-h',b_0),\ldots,(k-h,b_0),(a_1,b_1),\ldots,(a_m,b_m'),(k,k+1),(k-1,k+1),\ldots,(k-h',k+1)\]
in $\Gamma_{k}$. All conditions are checked as above.

{\bf Case 1.5.} This case represents the exception to Case 1.3, namely there is no $h'$ with $1\le h'\le h$ satisfying condition (\ref{case13}). This means that
\begin{equation}\label{case15}w(j+1)=w(j)+1\;\;\mbox{for}\;\;k-h\le j< k\,,\;\;\mbox{and}\;\;w(b_0)=w(k)+1\,.\end{equation}
The term on the right-hand side of (\ref{toprove}) corresponding to the chains $\gamma_i''$ and $\delta_i''$ in Case 1.3 is given by the chain 
\begin{equation}\label{ch15}(k,b_0),\ldots,(k-h,b_0),(a_1,b_1),\ldots,(a_m,b_m')\end{equation}
in $\Gamma_{k}$. All conditions are checked as above.

{\bf Case 1.6.} This case represents the exception to Case 1.4, namely condition (\ref{case15}) holds. The term on the right-hand side of (\ref{toprove}) corresponding to the chains $\gamma_i'''$ and $\delta_i'''$ in Case 1.4 is given by the chain 
\begin{equation}\label{ch16}(k,b_0),\ldots,(k-h,b_0),(a_1,b_1),\ldots,(a_m,b_m'),(k,k+1)\end{equation}
in $\Gamma_{k}$. All conditions are checked as above.

{\bf Summary of Case 1.} We have already accounted for most of the terms on the right-hand side of (\ref{toprove}). Let us specify precisely which ones in terms of the corresponding chains. We have accounted for all the chains in $\Gamma_{k-1}$ whose $k$-subchain is of length at least 1, cf. the definition of an $i$-subchain referring to the Pieri chain (\ref{pierichain}). We have also accounted for all the chains in $\Gamma_k$ whose $(k+1)$-subchain is of length at least 2. Finally, we have accounted for some chains in $\Gamma_k$ whose $(k+1)$-subchain is of length at most 1, namely the ones of the form (\ref{ch15}) and (\ref{ch16}). The latter requirement can be stated more concisely by saying that the initial subchain has length at least 2 and starts with $(k,\,\cdot\,)$. 

{\bf Case 2.}  This case corresponds to blocks of terms on the left-hand side of (\ref{toprove}) which cancel among themselves. We will use the chains $\gamma_i$ and $\delta_i$ introduced above. The crucial assumptions we make throughout this case, without mentioning them again, are that $m\ge 1$ and $a_1= k-h$.

Let us first consider the terms corresponding to the chains $\gamma_i$, for $i=1,\ldots,h+1$, and $\delta_i$, for $i=1,\ldots,h$, under the new assumptions. The computation of $m_p(\gamma_i)$ and $m_{p-1}(\gamma_i)$ is very similar to the one in Case 1.1; the only difference appears when $i=h+1$, when the last case in formula (\ref{coeff}) is used. To be more precise, we have
\begin{align*}\sum_{i=1}^{h+1}m_p(\gamma_i)&=(-1)^{m-p'-1}\binom{m-1}{p'}+(-1)^{m-p'}\binom{m}{p'}+(-1)^{m+1-p'}\binom{m}{p'-1}+\\&+\ldots+(-1)^{m+h-p'-2}\binom{m}{p'-h+2}+(-1)^{m+h-p'-1}\binom{m-1}{p'-h+1}\,,\\
\sum_{i=1}^{h+1}m_{p-1}(\gamma_i)&=(-1)^{m-p'}\binom{m-1}{p'-1}+(-1)^{m+1-p'}\binom{m}{p'-1}+\ldots+\\&\!\!\!\!\!\!\!\!\!\!\!\!\!\!\!\!\!\!\!\!\!\!\!+(-1)^{m+h-p'-2}\binom{m}{p'-h+2}+(-1)^{m+h-p'-1}\binom{m}{p'-h+1}+(-1)^{m+h-p'}\binom{m-1}{p'-h}.
\end{align*}
Hence, by using Pascal's identity, we have 
\[\sum_{i=1}^{h+1}m_p(\gamma_i)-\sum_{i=1}^{h+1}m_{p-1}(\gamma_i)=0\,.\]
Note that if $h=1$ then each of the two sums contains only the first and the last term in the corresponding formula. Therefore, both of them are 0 in this case, whereas only their difference is 0 when $h>1$. 

The calculation is completely similar for $\delta_i$, but now we must have $h\ge 2$. Finally, the same reasoning works for the class of chains obtained from $\gamma_i$ and $\delta_i$ by extending their Monk subchains via right concatenation with a fixed chain.

{\bf Case 3.}  Like Case 2, this one also corresponds to blocks of terms on the left-hand side of (\ref{toprove}) which cancel among themselves. The main difference is that, for the first time, we use the quantum edges of the quantum Bruhat graph. 

Let $\widetilde{\gamma}_i$ be the chain in $\widehat{\Gamma}_{k-1}$ of the form
\[(a_1,b_1),\ldots,(a_m,b_m'),(k-1,k),\ldots,(k-h,k)|(k-1,k),\ldots,(k-i,k)\,,\]
where $m\ge 0$, $h\ge 2$, $i=1,\ldots,h$, and (\ref{ab}) holds with $l=k-1$. Let $\widetilde{\delta}_i$ be the chain in $\widehat{\Gamma}_{k-2}$ of the form
\[(a_1,b_1),\ldots,(a_m,b_m'),(k-2,k-1),\ldots,(k-h,k-1)|(k-2,k),\ldots,(k-i,k)\,,\]
under the same conditions as above. It is easy to see that, for each $i$, the chains $\widetilde{\gamma}_i$ and $\widetilde{\delta}_i$ end in the same permutation (just commute the second transposition $(k-1,k)$ in $\widetilde{\gamma}_i$ to the left, in order to cancel the first one). 

{\bf Case 3.1.} Assume that $a_1<k-h$. We have
\begin{align*}
&m_p(\widetilde{\gamma}_i)=q_{k-1}m_{p-1}(\widetilde{\delta}_i)=(-1)^{m+h-p'-1}q_{k-1}\binom{m}{p'-h+1}\,,\\&m_{p-1}(\widetilde{\gamma}_i)=q_{k-1}m_{p-2}(\widetilde{\delta}_i)=(-1)^{m+h-p'}q_{k-1}\binom{m}{p'-h}\,.\end{align*}

{\bf Case 3.2.} Assume that we now have $a_1=k-h$. Then
\begin{align*}
&m_p(\widetilde{\gamma}_i)=q_{k-1}m_{p-1}(\widetilde{\delta}_i)=(-1)^{m+h-p'-1}q_{k-1}\binom{m-1}{p'-h+1}\,,\\&m_{p-1}(\widetilde{\gamma}_i)=q_{k-1}m_{p-2}(\widetilde{\delta}_i)=(-1)^{m+h-p'}q_{k-1}\binom{m-1}{p'-h}\,.\end{align*}

Hence, in both cases, each of the two terms on the left-hand side of (\ref{toprove}) corresponding to a chain $\widetilde{\gamma}_i$ cancels out with one of the two terms corresponding to the chain $\widetilde{\delta}_i$. Finally, note that the same reasoning works for the chains obtained from $\widetilde{\gamma}_i$ and $\widetilde{\delta}_i$ by extending their Monk subchains via right concatenation with a fixed chain (of length 1 or 2) starting with $(k,\,\cdot\,)$.

{\bf Case 4.} This case contains the exceptions to the previous ones. Consider a chain $\gamma$ of the form $(a_1,b_1),\ldots,(a_m,b_m')$, where $m\ge 1$ and (\ref{ab}) holds with modified upper bound for $a_1$, namely $a_1\le k-1$. Clearly $\gamma$, which has an empty Monk subchain, is in $\widehat{\Gamma}_{k-1}$ and, if $a_1\le k-2$, in $\widehat{\Gamma}_{k-2}$ as well. Let us also consider the chain
\[\widetilde{\gamma}:=\gamma,(k-1,k)|(k-1,k)\,.\]
Note that $\gamma$ coincides with the chain denoted by $\gamma_1$ in Case 1, asssuming that $h=0$, as well as with the chain $\widetilde{\delta}_1$ in Case 3, assuming that $h=1$. Furthermore, the chain $\widetilde{\gamma}$ coincides with the chain $\widetilde{\gamma}_1$ in Case 3, assuming that $h=1$.

{\bf Case 4.1.} Assume that $a_1<k-1$. Consider the terms on the left-hand side of (\ref{toprove}) corresponding to the chains $\gamma$ and $\widetilde{\gamma}$ viewed as chains in $\widehat{\Gamma}_{k-1}$ (two terms for each chain), as well as the two terms corresponding to the chain $\gamma$ viewed as a chain in $\widehat{\Gamma}_{k-2}$. We have
\begin{align*}&m_p(\gamma)+m_p(\widetilde{\gamma})=(-1)^{m-p'-1}\binom{m-1}{p'}+(-1)^{m-p'}q_{k-1}\binom{m}{p'}\,,\\
&m_{p-1}(\gamma)+m_{p-1}(\widetilde{\gamma})=(-1)^{m-p'}\binom{m-1}{p'-1}+(-1)^{m+1-p'}q_{k-1}\binom{m}{p'-1}\,,\\
&m_{p-1}(\gamma)=(-1)^{m-p'}\binom{m-1}{p'-1}\,,\;\;\;m_{p-2}(\gamma)=(-1)^{m+1-p'}\binom{m-1}{p'-2}\,;\end{align*}
here the four calculations correspond to the four sums on the left-hand side of (\ref{toprove}).  By Pascal's identity, the sum of the coefficients in the six terms corresponding to the chains $\gamma$ and $\widetilde{\gamma}$ is
\begin{equation}\label{sumcoeff}(-1)^{m-p'-1}\binom{m-1}{p'}-(-1)^{m-p'}\binom{m-1}{p'-1}=m_p(\gamma)-m_{p-1}(\gamma)\,.\end{equation}
But this is precisely the sum of the coefficients in the two terms on the right-hand side of (\ref{toprove}) corresponding to the chain $\gamma$.

The same reasoning can be repeated for the pairs of chains 
\begin{itemize}
\item $\gamma|(k,k+1)$ and $\widetilde{\gamma},(k,k+1)$;
\item $\gamma|(k,b_0)$ and $\widetilde{\gamma},(k,b_0)$;
\item $\gamma|(k,b_0),(k,k+1)$ and $\widetilde{\gamma},(k,b_0),(k,k+1)$.
\end{itemize}
Here we have $b_0>b_1$. The corresponding terms on the right-hand side of (\ref{toprove}) are given by the following chains in $\Gamma_k$, respectively:
\[\gamma,(k,k+1)\,;\;\;\;\;\;(k,b_0),\gamma\,;\;\;\;\;\;(k,b_0),\gamma,(k,k+1)\,.\]

{\bf Case 4.2.} Assume now that $a_1=k-1$. In this case, neither $\gamma$ nor $\widetilde{\gamma}$ belong to $\Gamma_{k-2}$. But we can again consider the terms on the left-hand side of (\ref{toprove}) corresponding to the chains $\gamma$ and $\widetilde{\gamma}$ viewed as chains in $\widehat{\Gamma}_{k-1}$ (two terms for each chain). We have
\begin{align*}&m_p(\gamma)+m_p(\widetilde{\gamma})=(-1)^{m-p'-1}\binom{m-1}{p'}+(-1)^{m-p'}q_{k-1}\binom{m-1}{p'}\,,\\
&m_{p-1}(\gamma)+m_{p-1}(\widetilde{\gamma})=(-1)^{m-p'}\binom{m-1}{p'-1}+(-1)^{m+1-p'}q_{k-1}\binom{m-1}{p'-1}\,;\end{align*}
here  the two calculations correspond to the first two sums on the left-hand side of (\ref{toprove}). It immediately follows that the sum of the coefficients in the four terms corresponding to the chains $\gamma$ and $\widetilde{\gamma}$ is precisely the one specified in (\ref{sumcoeff}). From this point, the reasoning is completely similar to the one in Case 4.1, including the discussion related to the three extra pairs of chains above. 

{\bf Summary of Case 4.} Let us identify the terms on the right-hand side of (\ref{toprove}) which were accounted for in Case 4. We will do this in terms of the corresponding chains. It is not hard to see that we accounted for all the chains in $\Gamma_{k-1}$ with empty $k$-subchain. The corresponding chains in $\Gamma_{k}$ are those characterized by: (i) their $(k+1)$-subchain has length at most 1; (ii) if they start with $(k,\,\cdot\,)$, then this transposition gives the whole initial subchain. The only difference between the chains in $\Gamma_{k-1}$ and $\Gamma_k$ treated in Cases 4.1 and 4.2 is that they do not or do contain a transposition $(k-1,\,\cdot\,)$, respectively. 

By comparing the chains on the right-hand side of (\ref{toprove}) accounted for in Cases 1 and 4 (see the Summary of those cases), it is easy to see that the two cases together cover all the chains and there are no overlaps. A careful analysis based on the cover condition also reveals that Cases 1-4 cover all the chains on the left-hand side of (\ref{toprove}), and again there are no overlaps. This completes the proof of Theorem \ref{main} in the case when $w$ is a Grassmannian permutation.

{\bf The general form of $w$.} If $w$ has one or more descents, all the chains considered above still start with a subchain of the form (\ref{pierichain}), but the second condition in (\ref{ab}) and (\ref{ab1}) do not necessarily hold. Instead, for any $i<j$, the following condition holds:
\[a_i'\le a_i\;\;\;\mbox{and either $a_j<a_i'$ (as before) or $a_j'\ge a_i$}\,.\]
Furthermore, in the Cases 1.3-1.6 above, the subchains of the chains on the left-hand side of (\ref{toprove}) consisting of transpositions of the form $(k,\cdot)$ can have lengths greater than 1 or 2; as a consequence, the corresponding chains $\gamma$ on the right-hand side of (\ref{toprove}) may have more than one subchain of the form $(k-h',b_0),\ldots,(k-h,b_0)$, and these subchains can be inserted in any position in the subchain $(a_1,b_1),\ldots,(a_m,b_m')$. These are the only differences from the Grassmannian case, and the reasoning is essentially the same as above.

We illustrate the case in which $w$ has more that one descent with an example. Let
\[w=1\,3\,5\,7\,9\,11\,12\,6\,10\,2\,8\,4\,,\;\;\;k=6\,,\;\;\;p=4\,.\]
Consider the chain
\[\overline{\gamma}:=(2,12),(1,10),(3,8),(2,8)\,.\]

{\bf Case 1.1.} We have
\[\begin{array}{lll}\gamma_1:=\overline{\gamma}|(5,6),(4,6)\;\;\;\;&\gamma_2:=\overline{\gamma},(5,6)|(4,6)\;\;\;\;&\gamma_3:=\overline{\gamma},(5,6),(4,6)|\\[0.05in]
\delta_1:=\overline{\gamma}|(5,6),(4,6)\;\;\;\;&\delta_2:=\overline{\gamma},(4,5)|(5,6)\,,&\end{array}\]
and 
\[{\rm end}(\gamma_i)={\rm end}(\delta_i)=2\,5\,6\,9\,11\,7\,12\,4\,10\,1\,8\,3.\]

Counting the appropriate markings of the above chains, we obtain
\[\begin{array}{llllll}\!\!m_4(\gamma_1)=0,&\!\!m_4(\gamma_2)=-1,&\!\!m_4(\gamma_3)=2,&\!\!m_3(\gamma_1)=-1,&\!\!m_3(\gamma_2)=2,&\!\!m_3(\gamma_3)=-1,\\[0.05in]
\!\!m_3(\delta_1)=-1,&\!\!m_3(\delta_2)=2,&&\!\!m_2(\delta_1)=1,&\!\!m_2(\delta_2)=-1\,.&\end{array}\]
For instance, in order to compute $m_4(\gamma_3)$, we start by observing that the (covers corresponding to the) transpositions $(2,12),(3,8),(5,6)$ must be marked (by the marking rules (M2) and (M3)), whereas $(2,8)$ must not be marked (by rule (M1)). The remaining transpositions $(1,10)$ and $(4,6)$ may or may not be marked, but we must mark precisely one of them in order to have a total of 4 markings. In each of the resulting two cases, the number of unmarked covers is even, namely 2. 

In consequence, we have
\begin{equation}\label{addcoeff}\frac{1}{1-q_5}\left(\sum_{i=1}^{3}m_4(\gamma_i)-\sum_{i=1}^{3}m_{3}(\gamma_i)-q_5\sum_{i=1}^{2}m_{3}(\delta_i)+q_5\sum_{i=1}^{2}m_{2}(\delta_i)\right)=1\,.\end{equation}
The chains $\gamma_i$ can be viewed as a single chain $\gamma$ in $\Gamma_{5}$. The coefficient in the corresponding term on the right-hand side of (\ref{toprove}), namely $-m_{3}(\gamma)$, is easily seen to also be 1.

{\bf Case 1.2.} In this case, we consider the chains $\gamma'$ and $\delta'$ obtained by appending to the chains $\gamma_i$ and $\delta_i$ the transposition $(6,7)$. By the same calculation as above, the sum of the coefficients in the terms on the left-hand side of (\ref{toprove}) corresponding to the chains $\gamma_i'$ and $\delta_i'$ is $-1$. The corresponding term on the right-hand side of (\ref{toprove}) is given by the chain $\gamma'=\overline{\gamma},(6,7),(5,7),(4,7)$ in $\Gamma_{6}$, and we have $m_4(\gamma')=-1$. 

{\bf Case 1.3.} In this case, we consider the chains $\gamma''$ and $\delta''$ obtained by appending to the chains $\gamma_i$ and $\delta_i$ the subchain $(6,11),(6,9)$. We have
\[{\rm end}(\gamma_i'')={\rm end}(\delta_j'')=2\,5\,6\,9\,11\,10\,12\,4\,8\,1\,7\,3.\]
The corresponding chain $\gamma''$ in $\Gamma_5$ on the right-hand side of (\ref{toprove}), which ends in the same permutation, is
\[\gamma''=(2,12),(4,11),(1,10),(5,9),(4,9),(3,8),(2,8),(5,6)\,.\]
The result of the calculation (\ref{addcoeff}) with $\gamma_i''$ and $\delta_i''$ instead of $\gamma_i$ and $\delta_i$ is identical, and we have $-m_3(\gamma'')=1$. 

{\bf Case 1.4.} is easily treated by combining Cases 1.2 and 1.3. The exceptions in Cases 1.5 and 1.6 are treated similarly.

{\bf Case 2.} In this case, we consider the chains 
\[\widehat{\gamma}_i=\gamma_i,(3,6),\;\,i=1,2,3,\;\;\;\;\;\widehat{\delta}_i=\delta_i,(3,6),\;\,i=1,2;\]
we also define $\widehat{\gamma}_4$ and $\widehat{\delta}_3$ in the obvious way. We have
\[m_j(\widehat{\gamma}_i)=m_j(\gamma_i),\;\,i=1,2,3,\;\,j=3,4\,,\]
as well as
\[m_4(\widehat{\gamma}_4)=-1,\;\;\;m_3(\widehat{\gamma}_4)=0\,.\]
Finally, we have
\[\sum_{i=1}^{4}m_4(\widehat{\gamma}_i)-\sum_{i=1}^{4}m_{3}(\widehat{\gamma}_i)=0-0=0\,.\]
The calculation for the chains $\widehat{\delta}_i$ is completely similar.

{\bf Case 3.1.} Let
\[\begin{array}{ll}\widetilde{\gamma}_1=\overline{\gamma},(5,6),(4,6)|(5,6)\;\;\;&\widetilde{\gamma}_2=\overline{\gamma},(5,6),(4,6)|(5,6),(4,6)\\[0.05in]
\widetilde{\delta}_1=\overline{\gamma},(4,5)|\;\;\;&\widetilde{\delta}_2=\overline{\gamma},(4,5)|(4,6)\end{array}\]
We have, for $i=1,2$,
\[m_4(\widetilde{\gamma}_i)=q_4 m_3(\widetilde{\delta}_i)=2q_4\,,\;\;\;\;m_3(\widetilde{\gamma}_i)=q_4 m_2(\widetilde{\delta}_i)=-q_4\,.\]
Hence, the terms corresponding to these chains in the left-hand side of (\ref{toprove}) cancel out.

An example in Case 3.2 is obtained by appending the transpositions $(3,6)$ and $(3,5)$ to the Pieri chains in $\widetilde{\gamma}_i$ and $\widetilde{\delta}_i$, $i=1,2$, respectively; we also have the extra chains $\widetilde{\gamma}_3$ and $\widetilde{\delta}_3$, but the calculations are completely similar. Finally, the exceptions in Case 4 are treated in a similar way to Cases 1 and 3.


\end{document}